\documentclass[draft]{amsart}

\DeclareMathOperator{\val}{val}
\DeclareMathOperator{\Hom}{Hom}

\DeclareMathOperator{\id}{id}
\DeclareMathOperator{\Ad}{Ad}
\DeclareMathOperator{\ad}{ad}

\DeclareMathOperator{\Aut}{Aut}
\DeclareMathOperator{\Spec}{Spec}
\DeclareMathOperator{\reg}{reg}
\DeclareMathOperator{\sep}{sep}
\DeclareMathOperator{\tame}{tame}
\DeclareMathOperator{\Gal}{Gal} 
\DeclareMathOperator{\cok}{cok}

\newcommand{\ox}{\mathcal O_X}

\usepackage{amscd}
\usepackage{amssymb}

\numberwithin{equation}{subsection}

\newtheorem{theorem}{Theorem}[subsection]
\newtheorem{corollary}[theorem]{Corollary}
\newtheorem{lemma}[theorem]{Lemma}
\newtheorem{proposition}[theorem]{Proposition}

\theoremstyle{definition}

\begin{document}
\title{Codimensions of root valuation strata}
\author{Mark Goresky}\thanks{The research of M.G. was 
supported in part by NSF grant
 DMS-0139986 and DARPA grant HR0011-04-1-0031}
\address{Goresky and MacPherson: School of Mathematics, Institute for
Advanced Study, Princeton, NJ, 08540 }
\author {Robert Kottwitz}\thanks{The research of R.K. was  
supported in part by NSF grants
 DMS-0071971 and DMS-0245639}
\address{Kottwitz: Department of Mathematics\\ University of Chicago\\ 5734
University Avenue\\ Chicago, Illinois 60637}
\author {Robert MacPherson}


\maketitle

\section{Introduction}
The topic of this paper arises naturally in the context of affine Springer
fibers, which we now take a moment to discuss. Let $G$ be a semisimple
complex algebraic group, and let $\mathfrak g$ denote its Lie algebra. We
then have the affine Grassmannian $X=G(F)/G(\mathcal O)$, where $\mathcal O$
is the ring $\mathbb C[[\epsilon]]$ of formal power series, and $F$ is its
fraction field $\mathbb C((\epsilon))$. 
For any $u \in \mathfrak g(F)=\mathfrak  g\otimes_\mathbb C F$ the closed
subset 
\[
X^u=\{g \in G(F)/G(\mathcal O): \Ad(g)^{-1}u \in \mathfrak g(\mathcal
O)=\mathfrak g \otimes_\mathbb C \mathcal O
\}
\]
of the affine Grassmannian, first studied by Kazhdan-Lusztig in
\cite{kazhdan-lusztig88}, is called the affine Springer fiber associated to
$u$. 

We now assume that $u$ is regular semisimple and write $T_u$ for its
centralizer in $G$, a maximal torus of $G$ over $F$. We will also need
$A_u$, the maximal $F$-split subtorus of $T_u$. If $u$ is integral, in the
sense that $\alpha(u)$ is integral over $\mathcal O$ for every root $\alpha$
of $T_u$, then $X^u$ is non-empty and may be viewed (see
\cite{kazhdan-lusztig88}) as the set of $\mathbb C$-points of a scheme
locally of finite type over $\mathbb C$. The dimension formula of
Bezrukavnikov-Kazhdan-Lusztig (see \cite{kazhdan-lusztig88} and
\cite{bezrukavnikov96}) states that 
\[
\dim X^u=(\delta_u -c_u)/2,
\]
where 
\begin{align*}
\delta_u&:=\val\det \bigl(\ad(u);\mathfrak g(F)/\mathfrak t_u(F)\bigr),\\
c_u&:=\dim T_u -\dim A_u.
\end{align*}
Here $\val$ is the usual valuation on $F$, normalized so that
$\val(\epsilon)=1$, and of course $\mathfrak t_u(F)$  denotes the Lie algebra
of the $F$-torus $T_u$. 

In particular $\dim X^u$ depends only on the  discrete invariant
$(\delta_u,c_u)$ of $u$. It is useful however to introduce a finer
invariant, still discrete in nature. For this we need to choose an algebraic
closure $\bar F$ of $F$. We denote by $\tau$ the unique element of
$\Gal(\bar F/F)$ that multiplies each $m$-th root of $\epsilon$ by $\exp(2\pi
i/m)$. Recall that $\tau$ is a topological generator of $\Gal(\bar F/F)$ and
allows us to identify that Galois group with the profinite completion of
$\mathbb Z$. 

Fix a maximal torus $T$ of $G$ over $\mathbb C$. We write $R$ for the set of
roots of $T$ in $G$, and $W$ for the Weyl group of $T$. Choose an element
$u' \in \mathfrak t(\bar F)$ that is $G(\bar F)$-conjugate to $u$. We attach
to $u'$ a pair $(w,r)$ in the following way: $w$ is the unique element of
$W$ such that $w\tau(u')=u'$, and $r:R \to \mathbb Q$ is the function
defined by $r(\alpha):=\val \alpha(u')$. Here we have extended our valuation
on $F$ to one on $\bar F$; the valuation of any $m$-th root of $\epsilon$
is then
$1/m$. Since $u$ is integral, the function $r$ takes values in the set of
non-negative rational numbers. The element $u'$ is not quite well-defined,
since it may be replaced by $xu'$ for any $x \in W$. This replaces $(w,r)$
by $(xwx^{-1},xr)$, where $(xr)(\alpha):=r(x^{-1}\alpha)$. 

All in all, we
have associated to $u$ a well-defined orbit $s$ of $W$ in the set of pairs
$(w,r)$, and $s$ is the desired discrete invariant of $u$. 
Clearly $s$ depends only on the $G(F)$-conjugacy class of $u$. Turning this
around, for a given orbit $s$, we let $\mathfrak g(F)_s$ denote the subset of
$\mathfrak g(F)$ consisting of all integral regular semisimple $u$ for which
the associated invariant is equal to $s$. 

Observe that the invariant $(\delta_u,c_u)$ can be expressed very simply 
in terms of the $W$-orbit of $(w,r)$. Indeed, we have 
\begin{align*}
\delta_u=\delta_r&:=\sum_{\alpha \in R} r(\alpha),\\ 
c_u=c_w&:=\dim \mathfrak t -\dim \mathfrak t^w,
\end{align*}
$\mathfrak t$ being the Lie algebra of $T$, and $\mathfrak t^w$ denoting the
fixed points of $w$ on $\mathfrak t$. Therefore the dimension of $X^u$ is
constant along each subset $\mathfrak g(F)_s$. 

We expect that something much stronger is true, namely that the cohomology
of $X^u$ is locally constant, in a suitable sense, along each subset
$\mathfrak g(F)_s$. In any case, this is true when the function $r$ is
constant (the equivalued case), as can be seen using the Hessenberg pavings
of \cite{gkm.pre2}.

Thus it is natural to study the subsets $\mathfrak g(F)_s$. This is best
done using the adjoint quotient $\mathbb A:=\mathfrak t/W$ and the natural
morphism 
\begin{equation}\label{eq.mor.in}
\mathfrak g(F) \to \mathbb A (F). 
\end{equation}
The set $\mathfrak g(F)_s$ is the preimage of a subset of $\mathbb
A(\mathcal O)$ that we will denote by $\mathbb A(\mathcal O)_s$. 

It is instructive to look at the case when $G=SL_2$.  Then $\mathbb A(F)=F$,
and the map \eqref{eq.mor.in} is 
\[
\det:\mathfrak{sl}_2(F) \to F.
\]
Each non-empty subset $\mathbb A(\mathcal O)_s$ is of the form
\[
Y_m=\{c \in \mathbb A(\mathcal O)=\mathcal O:\val c=m\}
\]
for some non-negative integer $m$. The pair $(w,r)$ corresponding to $m$ 
is determined as follows: $w$ is trivial (respectively, non-trivial) if $m$
is even (respectively, odd), and $r$ is the constant function with value
$m/2$. 

The subset $Y_m$ is admissible, in the sense that it is
the preimage of a subset in $\mathcal O/\epsilon^N \mathcal O$ once $N$ is
sufficiently large. This allows us to work with $Y_m$
just as if it were finite dimensional. In an obvious sense each $Y_m$ is 
(Zariski) locally closed, irreducible, non-singular
 of codimension $m$ in $\mathbb A(\mathcal O)$. 

One goal of this paper is to prove an analogous statement for any connected
reductive $G$ over an algebraically closed field $k$ in which the order of
the Weyl group is invertible. Theorem \ref{thm.a.str} says that
$\mathbb A(\mathcal O)_s$, when non-empty, is admissible, locally closed,
irreducible, and non-singular of codimension 
\[
d(w,r)+(\delta_r+c_w)/2
\] 
in $\mathbb A(\mathcal O)$. 
Here $\delta_r$, $c_w$ are the same integers as before, and $d(w,r)$ is the
codimension of $\mathfrak t_w(\mathcal O)_r$ in $\mathfrak t_w(\mathcal O)$,
where $\mathfrak t_w(\mathcal O)$ is the twist of $\mathfrak t$ by $w$, and
$\mathfrak t_w(\mathcal O)_r$ is a certain subset of $\mathfrak t_w(\mathcal
O)$ that maps onto $\mathbb A(\mathcal O)_s$ under 
\[
\mathfrak t_w(\mathcal O) \to \mathbb A(\mathcal O).
\]
The integer $d(w,r)$ is calculated in Proposition \ref{prop.t.str}(4). 

The second goal of the paper is to relate the geometry of $\mathfrak
t_w(\mathcal O)_r$ to that of $\mathbb A(\mathcal O)_s$ using the map 
\[
\mathfrak
t_w(\mathcal O)_r \twoheadrightarrow \mathbb A(\mathcal O)_s. 
\]
In Theorem \ref{thm.a.str} it is shown that $\mathfrak
t_w(\mathcal O)_r$ is smooth (in a suitable sense) over  $\mathbb A(\mathcal
O)_s$. Theorem \ref{thm.HH} gives a precise description of the
structure of this morphism. Combined   with Proposition
\ref{prop.t.str}, which concerns $\mathfrak
t_w(\mathcal O)_r$, it yields a clear picture of the structure of 
each individual stratum $\mathbb A(\mathcal O)_s$. 

However the methods of this paper shed little light on how the strata fit
together. We do not know, for example, whether the closure of $\mathbb
A(\mathcal O)_s$ is a union of strata. 

The paper contains some other results as well. We 
determine when $\mathbb
A(\mathcal O)_{(w,r)}$ is non-empty. Since $\mathfrak t_w(\mathcal O)_r$ 
maps
onto $\mathbb
A(\mathcal O)_{(w,r)}$, this is the same as determining when 
$\mathfrak t_w(\mathcal O)_r$ is non-empty, and this is done in Proposition
\ref{prop.ne.gen}. 

Now assume that 
$\mathbb A(\mathcal O)_{(w,r)}$ is non-empty. We show (Corollary
\ref{cor.wm=1}) that if 
$r$ takes values in
$\frac{1}{m}\mathbb Z$, then $w^m=1$. In particular, if $r$ takes 
values in $\mathbb Z$, then $w=1$. 
We also show (see subsection \ref{sub.eq.st}) that if the function $r$ is
constant, then the conjugacy class of $w$ is determined by $r$. (This is a
simple consequence of Springer's results \cite{springer74} on regular
elements in Weyl groups.) We do not know whether to expect that $w$ is
always redundant (more precisely, whether the non-emptiness of both $\mathbb
A(\mathcal O)_{(w,r)}$ and $\mathbb
A(\mathcal O)_{(w',r)}$ implies that $w$ and $w'$ are conjugate under some
element of the Weyl group that fixes $r$.) 

A substantial part of this work was done in June, 2000 at the  Centre
\'Emile Borel, which we would like to thank  both for its financial support
and the excellent working conditions it provided.   It is a pleasure to
thank M.~Sabitova for numerous helpful comments on a preliminary version of
this paper.

\section{Basic notation and definition of $\mathbb A(\mathcal O)'$} 
\subsection{Notation concerning $G$} Let $G$ be a connected reductive group
over an algebraically closed field $k$. We choose a maximal torus $T$ in
$G$, and write $\mathfrak t$ for its Lie algebra. 
Throughout this article we will assume that the order $|W|$ of the Weyl
group $W$ (of $T$ in $G$) is invertible in
$k$. 

We let $R \subset X^*(T)$ denote the set of roots of $T$ in $G$.
Occasionally we will need to fix a subset $R^+ \subset R$ of positive
roots. The differential of a root $\alpha$ is an element in the dual
space $\mathfrak t^*$ to $\mathfrak t$, and we will abuse notation a bit by
also writing $\alpha$ for this element of $\mathfrak t^*$. 

\subsection{Quotient variety $\mathbb A=\mathfrak t/W$} 
We will need the quotient variety $\mathbb A:=\mathfrak t/W$,  as well as
the canonical finite morphism 
\[
f:\mathfrak t \to \mathbb A.
\]
The notation $\mathbb A$ serves as a reminder that $\mathfrak t/W$ is
non-canonically isomorphic to affine $n$-space $\mathbb A^n$ with
$n=\dim(T)$. Indeed (see \cite{bourbaki.root}) the $k$-algebra of
$W$-invariant polynomial functions on $\mathfrak t$ is a polynomial algebra
on $n$ homogeneous generators $f_1,\dots,f_n$, called \emph{basic
invariants}. Choosing basic invariants $f_1,\dots,f_n$, we obtain a 
morphism 
\[
(f_1,\dots,f_n): \mathfrak t \to \mathbb A^n,
\]
which induces an isomorphism $\mathfrak t/W \cong \mathbb A^n$ and allows us
to view
$f$ as $(f_1,\dots,f_n)$. We will denote by $d_i$ the degree of the
polynomial $f_i$. 

\subsection{Open subsets of regular elements in $\mathfrak t$ and $\mathbb
A$} \label{sub.reg.jac}

Inside $\mathfrak t$ we have the $W$-invariant affine open subset
$\mathfrak t_{\reg}$ consisting of those elements $u \in \mathfrak t$ such
that $\alpha(u) \ne 0$ for all $\alpha \in R$. Since $|W|$ is invertible in
$k$, no root vanishes identically on $\mathfrak t$, and therefore
$\mathfrak t_{\reg}$ is non-empty. The quotient $\mathfrak t_{\reg}/W$ is a
non-empty affine open subset of $\mathbb A$ that we will denote by
$\mathbb A_{\reg}$. 

 Picking a basis in
the vector space $\mathfrak t$, we get coordinates $u_1,\dots,u_n$ on
$\mathfrak t$, and the Jacobian 
\[
J_u:=\det\Big( \frac {\partial f_i}{\partial u_j} \Bigr)
\]
is known (see \cite[Ch.~V, no.~5.5, Prop.~6]{bourbaki.root}) to have the
form 
\begin{equation}\label{eq.jac.pro}
J_u=c\prod_{\alpha \in R^+} \alpha(u)
\end{equation}
for some non-zero scalar $c \in k$.
In particular $\mathfrak t_{\reg}$ is  the set where the Jacobian
does not vanish, and therefore the restriction $f_{\reg}:\mathfrak t_{\reg}
\to \mathbb A_{\reg}$ of
$f$ is an \'etale covering with Galois group $W$. 

Later we will need the well-known
 identity \cite{bourbaki.root} 
\begin{equation}\label{eq.rt.d}
|R^+|=\sum_{i=1}^n(d_i-1),
\end{equation}
which can be proved by calculating the degree of the polynomial $J$ in two
different ways. 

\subsection{Definition of $\mathcal O$ and $F$} In fact we will mainly be
interested in $\mathbb A(\mathcal O)$, where $\mathcal O$ denotes the ring
$k[[\epsilon]]$ of formal power series. We also need the fraction field
$F=k((\epsilon))$ of $\mathcal O$. 

\subsection{Subset $\mathbb A(\mathcal O)'$ of $\mathbb A(\mathcal O)$} We
put $\mathbb A(\mathcal O)'=\mathbb A(\mathcal O) \cap \mathbb
A_{\reg}(F)$, the intersection being taken in $\mathbb A(F)$.  We stress
that this subset is considerably bigger than
$\mathbb A_{\reg}(\mathcal O)$. 
For example, when $G$ is $SL_2$, we have $\mathbb A(F)=F$, $\mathbb
A(\mathcal O)=\mathcal O$, $\mathbb A_{\reg}(F)=F^\times$, $\mathbb
A_{\reg}(\mathcal O)=\mathcal O^\times$, $\mathbb A(\mathcal O)'=\mathcal O
\setminus \{0\}$.  Our first task in this paper is to partition the set
$\mathbb A(\mathcal O)'$. Roughly speaking, this involves  two ingredients: 
valuations of roots and Weyl group elements. We begin by discussing
valuations of roots.

\section{Valuations of roots: split case}

\subsection{Normalization of the valuation on $F$} 
We normalize the valuation on $F$  so that
$\val(\epsilon)=1$. 

\subsection{Definition of $\mathfrak t (\mathcal O)'$} 
Put  $\mathfrak
t(\mathcal O)':=\mathfrak t(\mathcal O) \cap \mathfrak t_{\reg}(F)$.

\subsection{Definition of $r_u$}\label{sub.def.ru} 
For any $u \in \mathfrak
t (\mathcal O)'$ we define a function $r_u$ on $R$ by 
\[
r_u(\alpha)=\val \alpha(u)
\] 
for each root $\alpha$. 
It is clear that $r_u$ takes values in the set of non-negative integers. 

Since $W$ acts on $R$, it acts on functions $r$ on $R$ by the rule
$(wr)(\alpha)=r(w^{-1}\alpha)$. It is clear that 
\begin{equation}\label{eq.r.w.u}
r_{wu}=wr_u
\end{equation}
for all $w \in W$ and $u \in \mathfrak t(\mathcal O)'$. 

\subsection{Properties of the function $r_u$} \label{sub.r_u}
 Let $u \in \mathfrak t
(\mathcal O)'$. It is obvious that
\begin{equation}\label{eq.-sym}
r_u(-\alpha)=r_u(\alpha).
\end{equation}
However the non-archimedean property of valuations gives much more than
this, as we will now see. 

Fix some function $r$ on $R$ with values in the set of non-negative
integers. We define a subset $\mathfrak t(\mathcal O)_r$ of $\mathfrak
t(\mathcal O)'$ by 
\[
\mathfrak t(\mathcal O)_r:=\{u \in \mathfrak
t(\mathcal O)' : r_u=r \}.
\]
We also use 
$r$ to define a chain 
\[
R=R_0 \supset R_1 \supset R_2 \supset R_3 \supset \dots
\] 
of subsets    
\[
R_m:=\{\alpha \in R: r(\alpha) \ge m \}.
\] 
We will need  the linear subspaces
\[
\mathfrak a_m:=\{u \in \mathfrak t:\alpha(u)=0 \quad \forall \, \alpha \in
R_m \}.
\]
These form an increasing chain 
\[
\mathfrak a_0 \subset \mathfrak a_1 \subset \mathfrak a_2 \subset \dots
\]
with $\mathfrak a_m=\mathfrak t$ for large enough $m$. Finally, for each
$m \ge1$ we will need the subset 
\[
\mathfrak a_m^\sharp:=\{u \in \mathfrak a_m:\alpha(u)\ne0\quad\forall\,
\alpha \in R_{m-1}\setminus R_{m}\}
\] 
of $\mathfrak a_m$.

\begin{proposition}\label{prop.split.empty}
The set $\mathfrak t(\mathcal O)_r$ is non-empty if and only
if each subset $R_m$ is $\mathbb Q$-closed, in the sense that if $\alpha
\in R$ is a $\mathbb Q$-linear combination of elements in $R_m$, then
$\alpha$ itself lies in $R_m$. Moreover $\mathfrak t(\mathcal O)_r$ has the
following description:  
$u \in \mathfrak t(\mathcal O)$ lies in $\mathfrak t(\mathcal O)_r$
if and only if the coefficients $u_j$ in the power series expansion of $u$ 
satisfy 
$
u_j \in \mathfrak a_{j+1}^\sharp
$ 
for all $j \ge 0$. 
\end{proposition}
\begin{proof} $(\Longrightarrow)$ Choose $u \in \mathfrak t(\mathcal O)_r$
and  expand it as a formal power series 
\[
u=\sum_{j=0}^\infty u_j \epsilon^j
\]
with coefficients $u_j \in \mathfrak t$.  
Clearly 
$
R_m=\{\alpha \in R: \alpha (u_i)=0
\quad\forall \, i=0,\dots,m-1\}. 
$ 
 It now follows from Proposition \ref{prop.stein} that $R_m$ is $\mathbb
Q$-closed. 

$(\Longleftarrow)$ Assuming that each $R_m$ is $\mathbb
Q$-closed,  we must show that $\mathfrak t(\mathcal O)_r$ is non-empty.
It is clear from the definitions that an
element $u \in \mathfrak t(\mathcal O)$ lies in $\mathfrak t(\mathcal O)_r$
if and only if the coefficients $u_j$ in its power series expansion satisfy 
$
u_j \in \mathfrak a_{j+1}^\sharp.
$
Thus we just need to show that each $\mathfrak a_{j+1}^\sharp$ is
non-empty. Since $R_{j+1}$ is $\mathbb Q$-closed, it is the root system
$R_M$ of some Levi subgroup $M \supset T$ (see the proof of  
Proposition
\ref{prop.stein}(3)).  Lemma
\ref{lem.goth.am} then tells us  that  no root
in
$R_j\setminus R_{j+1}$ vanishes identically on $\mathfrak a_{j+1}$, from
which it follows immediately that $\mathfrak a^\sharp_{j+1}$ is
non-empty.   
\end{proof}

\section{Twisted forms $\mathfrak t_w(\mathcal O)$  and strata $\mathfrak
t_w(\mathcal O)_r$} The subsets $\mathfrak t(\mathcal O)_r$ will help us to
understand $\mathbb A(\mathcal O)'$, but they are not enough, since the
canonical map $\mathfrak t(\mathcal O)' \to \mathbb A(\mathcal O)'$ is by
no means surjective. 
 In order to get a handle on all elements of
$\mathbb A(\mathcal O)'$ we need some twisted forms $\mathfrak t_w(\mathcal
O)$ of
$\mathfrak t$ over
$\mathcal O$. 

 For example, when $G$ is $SL_2$ (so that $2$ is required to be invertible
in $k$), the map in question  is---up to
multiplication by a scalar in $k^\times$---the squaring map from $\mathcal
O\setminus \{0\}$ to $\mathcal
O\setminus \{0\}$, whose image consists precisely of those elements in
$\mathcal O$ with even valuation.  
To obtain the missing elements  we need to replace $\mathfrak
t(\mathcal O)=\mathcal O$ by the $\mathcal O$-module of elements in
$k[[\epsilon^{1/2}]]$ having trace $0$ in $k[[\epsilon]]$, or, in other
words, the 
$\mathcal O$-module (free of rank $1$) $\mathcal O\epsilon^{1/2}$. The
squares of the non-zero elements in $\mathcal O\epsilon^{1/2}$ then yield
all elements in $\mathcal O$ having odd valuation. The $\mathcal O$-module
$\mathcal O\epsilon^{1/2}$ will turn out to be the twisted form $\mathfrak
t_w(\mathcal O)$ obtained from the non-trivial element  $w \in W$.

We begin by  reviewing tamely ramified extensions
of $F$. Next we  define $\mathfrak t_w(\mathcal O)$. Then we use
valuations of roots to define subsets 
$\mathfrak t_w(\mathcal O)_r$ of $\mathfrak t_w(\mathcal O)$.  Finally we
determine when the strata $\mathfrak t_w(\mathcal O)_r$ are non-empty.

\subsection{Review of $F_{\tame}$} 

We now need to choose an algebraic
closure $\bar F$ of $F$. We denote by $F_{\sep}$ the separable closure  of
$F$ in $\bar F$, and by $F_{\tame}$ the maximal tamely ramified extension 
of $F$ in $F_{\sep}$.  

It is well-known that $F_{\tame}$ has the following concrete description.
For any positive integer $l$ that is invertible in $k$, we choose an
$l$-th root $\epsilon^{1/l}$ of $\epsilon$ in $\bar F$, and we do this in
such a way that $(\epsilon^{1/lm})^m=\epsilon^{1/l}$ for any two positive
integers $l,m$ that are both invertible in
$k$. The field
$F_l:=F(\epsilon^{1/l})=k((\epsilon^{1/l}))$ is cyclic of degree $l$ over
$F$, and is independent of the choice of $l$-th root of $\epsilon$. Moreover
$F_{\tame}$ is the union of all the subfields $F_l$.  

For any positive integer $l$ that is invertible in $k$, we also choose a
primitive 
$l$-th root $\zeta_l$ of $1$ in $k$, and  we do this in such a
way that $(\zeta_{lm})^m=\zeta_l$ for any two positive integers $l,m$ that
are both invertible in
$k$. We use $\zeta_l$ to obtain a generator $\tau_l$
of 
$\Gal(F_l/F)$, namely the unique automorphism of $F_l/F$ taking 
$\epsilon^{1/l}$ to $\zeta_l \epsilon^{1/l}$. These generators are
consistent with each other as $l$ varies, and therefore fit together to
give an automorphism $\tau_\infty$ of $F_{\tame}/F$ whose restriction to 
each 
$F_l$ is $\tau_l$. Clearly $\tau_\infty$ is a topological generator of the
topologically cyclic group $\Gal(F_{\tame}/F)$.

\subsection{Definition of $\mathfrak t_w(\mathcal O)$}
Now we can construct the twisted forms of $\mathfrak t$ alluded to before. 
To get such a twist we need to start with an element $w \in W$. We then
take $l$ to be the order $o(w)$ of $w$, a positive integer that is
invertible in $k$. We write $E$ instead of $F_l$ and $\tau_E$ instead of
$\tau_l$. Moreover we write $\epsilon_E$ for $\epsilon^{1/l}$, so that
$E=k((\epsilon_E))$ and the valuation ring $\mathcal O_E$ in $E$ is
$k[[\epsilon_E ]]$. 

Then we put 
\begin{equation}\label{eq.def.twO}
\mathfrak t_w(\mathcal O):=\{u \in \mathfrak t(\mathcal O_E) : w\tau_E(u)=u
\}.
\end{equation}
More generally, for any $\mathcal O$-algebra $A$,  we put 
\[
\mathfrak t_w(A):=\mathfrak t_w(\mathcal O) \otimes_\mathcal O A.
\]
Since it will become clear in subsection \ref{sub.tsO}  that
$\mathfrak t_w(\mathcal O)$ is a free $\mathcal O$-module of rank $n$, where
$n=\dim_k\mathfrak t$, we see that 
$\mathfrak t_w$ is a scheme over
$\mathcal O$  isomorphic to affine
$n$-space over $\mathcal O$. 

Note that only the conjugacy class of $w$ in $W$ really matters: given $x
\in W$ we obtain an isomorphism $u \mapsto xu$ from $\mathfrak t_w(\mathcal
O)$ to $\mathfrak t_{xwx^{-1}}(\mathcal
O)$. This shows too that the centralizer $W_w$ (of $w$ in $W$) acts on
$\mathfrak t_w(\mathcal O)$ (and hence on $\mathfrak t_w$ over $\mathcal
O$). 

\subsection{Description of $\mathfrak t_w(\mathcal O)$}\label{sub.tsO}
It is easy to describe $\mathfrak t_w(\mathcal O)$ in terms of the
eigenspaces for the action of $w$ on
$\mathfrak t$. Since $w$ has order $l$,  the only possible
eigenvalues are
$l$-th roots of unity. Because $l$ is invertible in $k$, we then have 
\[
\mathfrak t=\bigoplus_{j=0}^{l-1} \mathfrak t(w,j),
\]
where $\mathfrak t(w,j)$ denotes the eigenspace 
$\mathfrak t(w,j):=\{v \in \mathfrak t:wv=\zeta_l^{-j}v\}$.

An element $u \in \mathfrak t(\mathcal O_E)$ can be expanded as a formal
power series 
\[
\sum_{j=0}^\infty u_j \epsilon_E^j
\]
with $u_j \in \mathfrak t$,  and we see from \eqref{eq.def.twO} that $u \in
\mathfrak t_w(\mathcal O)$ if and only if  $u_j \in \mathfrak
t(w,j)$ for all $j \ge 0$.  
Thus there is a canonical $\mathcal O$-module isomorphism   
\begin{equation}\label{eq.tw.eig}
\mathfrak t_w(\mathcal O)\cong \bigoplus_{j=0}^{l-1} \mathcal O\epsilon^j_E
\otimes_k
\mathfrak t(w,j).
\end{equation}
%Under this isomorphism the element $1\otimes u_j$ (with 
% $u_j \in \mathfrak t(w,j)$ and $0 \le j \le l-1$) of the right side
%corresponds to the element $u_j\epsilon_E^j$ of the left side. 

\subsection{Description of $\mathfrak t_w$ as a fixed point scheme}
\label{sub.tw.fp}
We write $R_{\mathcal O_E/\mathcal O} \mathfrak t$ for the scheme over
$\mathcal O$ obtained by starting with $\mathfrak t$, then extending scalars
from $k$ to $\mathcal O_E$, then (Weil) restricting scalars from $\mathcal
O_E$ to $\mathcal O$. For any $\mathcal O$-algebra $A$ we then have 
\[
(R_{\mathcal O_E/\mathcal O} \mathfrak t)(A)=\mathfrak t(A \otimes_\mathcal
O \mathcal O_E).
\]
Of course $R_{\mathcal O_E/\mathcal O} \mathfrak t$ is non-canonically
isomorphic to affine space of dimension ${nl}$ over $\mathcal O$. 

The automorphism $\tau_E$ of $\mathcal O_E/\mathcal O$ induces an
automorphism 
\[
\tau_E: R_{\mathcal O_E/\mathcal O} \mathfrak t \to R_{\mathcal
O_E/\mathcal O} \mathfrak t
\]
(given on $A$-valued points  by the map induced by the $\mathcal O$-algebra
automorphism $\id_A \otimes\tau_E$ of $A \otimes_\mathcal O \mathcal O_E$).
Moreover our $W$-action on $\mathfrak t$ induces a $W$-action on
$R_{\mathcal O_E/\mathcal O} \mathfrak t$. The actions of $W$ and $\tau_E$
 commute, and therefore the cyclic
group $\mathbb Z/l\mathbb Z$ acts on $R_{\mathcal O_E/\mathcal O} \mathfrak
t$ with the standard generator of that cyclic group acting by $w\circ
\tau_E$. 

Using \eqref{eq.tw.eig}, one sees easily that for 
 any $\mathcal O$-algebra $A$ we have 
\begin{equation}\label{eq.TWA}
\mathfrak t_w(A)=\{u \in \mathfrak t(A \otimes_\mathcal O \mathcal O_E)  
  : w\tau_E(u)=u \},
\end{equation} 
and hence that $\mathfrak t_w$ is the fixed point scheme (see appendix
\ref{sec.fp.gx}) of the action of $\mathbb Z/l\mathbb Z$ on 
$R_{\mathcal O_E/\mathcal O} \mathfrak t$. As
a special case  of \eqref{eq.TWA} we have
\[
\mathfrak t_w(F)=\{ u \in \mathfrak t(E): w\tau_E(u)=u \}.
\]

\subsection{Definition of $\mathfrak t_w(\mathcal O)'$} We put 
\[
\mathfrak t_w(\mathcal O)':=\mathfrak t(\mathcal O_E)' \cap 
\mathfrak t_w(\mathcal O).  
\]
Thus $u \in \mathfrak t_w(\mathcal O)$ lies in $\mathfrak t_w(\mathcal
O)'$  if and only if $\alpha(u)\ne 0$ for all $\alpha \in R$.

\subsection{Definition of strata $\mathfrak t_w(\mathcal O)_r$ in $\mathfrak
t_w(\mathcal O)'$} \label{sub.ne.str}
We extend the valuation on the field $F$ to a valuation, still denoted
$\val$, on
$\bar F$. In particular we have $\val(\epsilon^{1/l})=1/l$. 

Let $\mathcal R$ denote the set of functions on $R$ with values in the set
of non-negative rational numbers.  For $r \in \mathcal R$ we put 
\[
\mathfrak t_w(\mathcal O)_r:=\{ u \in \mathfrak t_w(\mathcal O): \val
\alpha(u)=r(\alpha) \quad \forall \, \alpha \in R \}. 
\]
It is clear that $\mathfrak t_w(\mathcal O)'$ is the disjoint union of the
strata $\mathfrak t_w(\mathcal O)_r$, many of which are empty. 

The Weyl group acts on itself by conjugation, and it also acts on $\mathcal
R$ (see subsection \ref{sub.def.ru}); thus we have an action of $W$ on the
set of pairs
$(w,r) \in W \times \mathcal R$.  Note that only the
$W$-orbit  of
$(w,r)$  really matters: given
$x
\in W$ we obtain an isomorphism $u \mapsto xu$ from $\mathfrak t_w(\mathcal
O)_r$ to $\mathfrak t_{xwx^{-1}}(\mathcal
O)_{xr}$. 

\subsection{Freeness of the $W_w$-action on $\mathfrak t_w(\mathcal
O/\epsilon^N\mathcal O)_{r<N}$} \label{sub.free} 
The centralizer $W_w$ acts freely on $\mathfrak t_w(\mathcal O)'$
by Proposition \ref{prop.stein}. Now let $N$ be a positive integer. We are
going to define an open  subset $\mathfrak t_w(\mathcal
O/\epsilon^N\mathcal O)_{r<N}$ of the $k$-variety $\mathfrak t_w(\mathcal
O/\epsilon^N\mathcal O)$ (see \ref{sub.gre.sm}) on which $W_w$ acts freely.
Here is the definition:
\[
\mathfrak t_w(\mathcal
O/\epsilon^N\mathcal O)_{r<N}:=\{u \in \mathfrak t_w(\mathcal
O/\epsilon^N\mathcal O): \alpha(u)\ne 0 \quad \forall \, \alpha \in R\}.
\]
(Note that $\alpha(u)$ is an element of the ring $\mathcal
O_E/\epsilon^N\mathcal O_E$.) The set   $\mathfrak t_w(\mathcal
O/\epsilon^N\mathcal O)_{r<N}$ can also be described as the image in $\mathfrak t_w(\mathcal
O/\epsilon^N\mathcal O)$ of all strata $\mathfrak t_w(\mathcal O)_r$ for
which $r$ satisfies the condition  $r(\alpha) < N$ for all $\alpha \in
R$. 

Now we verify that  $W_{w}$ acts freely on $\mathfrak t_w(\mathcal
O/\epsilon^N\mathcal O)_{r<N}$. Let $ u \in \mathfrak t_w(\mathcal
O/\epsilon^N\mathcal O)_{r<N}$ and
expand it as 
\[
 u=\sum_{j=0}^{Nl-1} u_j \epsilon_E^j. 
\]
Suppose that some element $x \in W_{w}$ fixes $u$. Then $x$ fixes each
coefficient $u_j$. It follows from Proposition \ref{prop.stein} (1) that
$x$ lies in the Weyl group of the root system consisting of all roots
$\alpha \in R$ such that $\alpha(u_j)=0$ for all $j$. Since $\alpha(u)\ne 0
$ for all $\alpha \in R$, there are no such roots, and therefore $x=1$. 

\subsection{Which strata are non-empty?}\label{sub.wh.ne}
 We are now going to determine
which strata $\mathfrak t_w(\mathcal O)_r$ are non-empty. (Only the
$W$-orbit of $(w,r)$ matters.) We begin by listing some useful necessary
conditions. 

Let $u \in \mathfrak
t_w(\mathcal O)'$. Then $\alpha(u) \in \mathcal O_E$ and thus 
$\val\alpha(u) \in \frac{1}{l}\mathbb Z$. Therefore a necessary
condition for non-emptiness of $\mathfrak t_w(\mathcal O)_r$ is that 
$r$ take values in
$\frac{1}{l}\mathbb Z$. Of course this statement can be sharpened a little,
since a particular root
$\alpha$ may be defined over $k((\epsilon^{1/l'}))$ for some divisor $l'$
of $l$, in which case it is necessary that $r(\alpha)$  lie in 
$\frac{1}{l'}\mathbb Z$. 

Now assume that $r$ does take values in $\frac{1}{l}\mathbb Z$, and define
an integer valued function $r_E$ on $R$ by $r_E(\alpha)=lr(\alpha)$. It is
clear from the definitions that 
\[
\mathfrak t_w(\mathcal O)_r=\mathfrak t(\mathcal O_E)_{r_E} \cap
\mathfrak t_w(\mathcal O).
\]
(Here we are applying definitions we have already made for $F$ to the field
$E$, so that when interpreting the right side of this equality one should
be thinking of the normalized valuation on $E$, rather than the one that
extends the valuation on $F$. That is why we need $r_E$ instead of $r$.) Now
Proposition
\ref{prop.split.empty} tells us exactly when
$\mathfrak t(\mathcal O_E)_{r_E}$ is non-empty. We conclude that another
necessary condition for the non-emptiness of $\mathfrak t_w(\mathcal O)_r$
is that the subset 
\[
R_m:=\{\alpha \in R: r(\alpha) \ge m/l \}
\]
of $R$ be $\mathbb Q$-closed for every non-negative integer $m$. We now
assume that this condition on $r$ also holds. 

Our stratum $\mathfrak t_w(\mathcal O)_r$ might still be empty. To settle
the question we need once again to consider  
the  vector spaces 
\[
\mathfrak a_m:=\{u \in \mathfrak t:\alpha(u)=0 \quad \forall \, \alpha \in
R_m \}
\]
 and their open subsets 
\[
\mathfrak a_m^\sharp:=\{u \in \mathfrak a_m:\alpha(u)\ne0\quad\forall\,
\alpha \in R_{m-1}\setminus R_{m}\}.
\] 
\begin{lemma}\label{lem.t.str.ne}
Let $u \in\mathfrak t(\mathcal O_E)$, and expand $u$ as
a power series  
\[
\sum_{j=0}^\infty u_j \epsilon_E^j
\]
with $u_j \in \mathfrak t$.  
Then $u \in \mathfrak t_w(\mathcal O)_r$ if and only if $u_j \in \mathfrak
t(w,j) \cap \mathfrak a_{j+1}^\sharp$ for all $j \ge 0$. Consequently
$\mathfrak t_w(\mathcal O)_r$ is non-empty if and only if $\mathfrak
t(w,j) \cap \mathfrak a_{j+1}^\sharp$ is non-empty for all $j \ge 0$. 
\end{lemma}

\begin{proof}
We  observed in \ref{sub.tsO} that $u \in
\mathfrak t_w(\mathcal O)$ if and only if  $u_j \in \mathfrak
t(w,j)$ for all $j \ge 0$. Moreover it follows  from Proposition
\ref{prop.split.empty} that $u \in \mathfrak t(\mathcal O_E)_{r_E}$ if and
only if $u_j \in \mathfrak a_{j+1}^\sharp$ for all $j \ge 0$. 
\end{proof}

We can reformulate the non-emptiness result in the last lemma in a
slightly better way, but for this we first need to note that there is
another obvious necessary condition for non-emptiness. Indeed, for any $u
\in \mathfrak t(\mathcal O_E)'$ we have $r_{\tau_E(u)}=r_u$ (obvious) and
hence $r_{w\tau_E(u)}=wr_u$ (use \eqref{eq.r.w.u}). It follows that if $u
\in \mathfrak t_w(\mathcal O)'$, then $r_u=wr_u$. 

Thus we see that if $\mathfrak t_w(\mathcal O)_r$ is non-empty, then $w$
stabilizes $r$, from which it follows that $w$ stabilizes the subsets $R_m$
and the subspaces $\mathfrak a_j$ of $\mathfrak t$, so that we obtain an
action of $w$ on each quotient space $\mathfrak a_{j+1}/\mathfrak a_j$. In
this situation we may consider the eigenspace 
\[
(\mathfrak
a_{j+1}/\mathfrak a_j)(w,j):=\{v \in \mathfrak a_{j+1}/\mathfrak a_j :
wv=\zeta_l^{-j}v \}.
\]
Let us also note that since each root in $R_{j}$ vanishes identically on
$\mathfrak a_j$, our subset $\mathfrak a_{j+1}^\sharp$ is the preimage
under $\mathfrak a_{j+1} \to \mathfrak a_{j+1}/\mathfrak a_j$ of the set 
\[
(\mathfrak a_{j+1}/\mathfrak a_j)^\sharp:=
\{u \in \mathfrak a_{j+1}/\mathfrak a_j:\alpha(u)\ne0\quad\forall\,
\alpha \in R_{j}\setminus R_{j+1}\}.
\] 
\begin{proposition}\label{prop.ne.gen} 
The stratum 
$\mathfrak t_w(\mathcal O)_r$ is non-empty if and only if 
the following four  conditions hold.
\begin{enumerate}
\item $r$ takes values in $\frac{1}{l}\mathbb Z$.
\item $R_m$ is $\mathbb Q$-closed for all $m \ge 0$. 
\item $wr=r$.
\item  $(\mathfrak
a_{j+1}/\mathfrak a_j)(w,j) \cap (\mathfrak a_{j+1}/\mathfrak a_j)^\sharp$
is non-empty  for all $j \ge 0$.
\end{enumerate}
 \end{proposition}

\begin{proof}
This follows from the previous lemma, since $\mathfrak t(w,j) \cap \mathfrak
a_{j+1}$ projects onto $(\mathfrak
a_{j+1}/\mathfrak a_j)(w,j)$. 
\end{proof}

\begin{proposition}
Suppose that $\mathfrak t_w(\mathcal O)_r$ is non-empty, and suppose that
$r(\alpha) \in \mathbb Z$ for all $\alpha \in R$. Then $w=1$. 
\end{proposition}
\begin{proof}
Let $u \in \mathfrak t_w(\mathcal O)_r$ and expand $u$ as a power series 
\[
\sum_{j=0}^\infty u_j \epsilon_E^j. 
\]
We are going to apply Proposition
\ref{prop.stein} to the subset $S:=\{u_{ml}:m=0,1,2,\dots\}$. Since
$u_{ml} \in \mathfrak t(w,ml)$, we see that $w$ fixes each element of $S$
and therefore lies in the subgroup $W_S=W(R_S)$ of Proposition
\ref{prop.stein}. However $R_S$ is empty, since for any $\alpha \in R$ we
have $\alpha(u_{r(\alpha)l})\ne 0$. Therefore $w=1$.   
\end{proof}

\begin{corollary}\label{cor.wm=1} 
Let $m$ be a positive integer. 
Suppose that $\mathfrak t_w(\mathcal O)_r$ is non-empty, and suppose that
$r(\alpha) \in \frac{1}{m}\mathbb Z$ for all $\alpha \in R$. Then $w^m=1$.
\end{corollary}
\begin{proof}
The idea is to extend scalars from $F$ to $F_m$. We denote the valuation
ring in $F_m$ by $\mathcal O_m$. Then we have 
\[
 \mathfrak t_w(\mathcal O)_r \subset \mathfrak t_{w^m}(\mathcal
O_m)_{mr},
\]
which shows that $\mathfrak t_{w^m}(\mathcal
O_m)_{mr}$ is non-empty. Since $mr$ takes integral values, the previous
result, applied to $F_m$ rather than $F$, tells us that
$w^m=1$. 
\end{proof}

\subsection{Equivalued strata}\label{sub.eq.st}
 We say that a stratum
$\mathfrak t_w(\mathcal O)_r$ is \emph{equivalued} if the function $r$ is
constant.  We will now use Proposition \ref{prop.ne.gen} to reduce the
problem of classifying non-empty equivalued strata to a problem that has
already been solved by Springer \cite{springer74}. Because Springer works
over a base-field of characteristic $0$, we temporarily do so too, just in
this subsection. (It seems quite likely that our usual hypothesis that
$|W|$ be invertible in $k$ suffices, but we have not checked this
carefully.) 

We need to recall the following
definition (due to Springer): an element $w \in W$ is said to be
\emph{regular} if there exists a non-zero eigenvector $u$ of $w$ in
$\mathfrak t$ that is regular (in the sense that no root $\alpha$ vanishes
on $u$). 
When $w$ is regular of order $l$,
Springer
\cite{springer74} shows that the eigenspace $\mathfrak t(w,j)$ contains a
regular element of $\mathfrak t$ if and only if $j$ is relatively prime to
$l$. 

\begin{proposition}
Let $a/b$ be a non-negative rational number written in least common
terms, so that $b$ is positive and $(a,b)=1$. Let $r$ be the constant
function on $R$ with value $a/b$, and let $w \in W$.  Then $\mathfrak
t_w(\mathcal O)_r$ is non-empty if and only if $w$ is  regular of order
$b$. 
\end{proposition}
\begin{proof} We use our usual notation. In particular $l$ denotes the
order of $w$. From Proposition \ref{prop.ne.gen} we see that $\mathfrak
t_w(\mathcal O)_r$ is non-empty if and only $b$ divides $l$ and
the eigenspace $\mathfrak t(w,al/b)$ contains a regular element of
$\mathfrak t$. The proposition follows from this and the result of Springer
mentioned just before the statement of the proposition. 
\end{proof}

Let us now recall a beautiful result from Springer's paper (see
\cite[Theorem 4.2]{springer74}, as well as the remarks following the
proof of that theorem):  all regular elements in $W$ of a given order are
conjugate. Combining this with the previous proposition, we see that for
a given $a/b$ in least common terms (letting $r$ denote, as before, the
constant function with value $a/b$), there are two possibilities. The first
is that there is no regular element of $W$ having order $b$. In this case
there are no non-empty strata $\mathfrak t_w(\mathcal O)_r$. The second is
that there are 
 regular elements in $W$ having order $b$, in which case there is a single
conjugacy class of such elements, and $\mathfrak t_w(\mathcal O)_r$ is
non-empty if and only if $w$ lies in this conjugacy class. The somewhat
surprising conclusion is that $w$ is essentially redundant: given $a/b$,
there is at most one $W$-orbit of pairs $(w,r)$ for which $r$ is the
constant function with value $a/b$ and $\mathfrak t_w(\mathcal O)_r$ is
non-empty. 

This raises an obvious question. Let $r$ be a non-negative rational valued
function on
$R$, and let $W_r$ denote the stabilizer of $r$ in $W$. Suppose that 
$\mathfrak t_w(\mathcal O)_r$ and $\mathfrak t_{w'}(\mathcal O)_r$ are both
non-empty. Is it then true that $w$, $w'$ are conjugate under $W_r$? This
is a question, not a conjecture. We do not know whether to expect a
positive or negative answer. 

We should also remark that  Springer \cite{springer74} gives a list of
the regular elements in the Weyl group of each irreducible root system.
Together with the proposition we just proved, this gives a classification
of all  non-empty equivalued strata. Take $G$ to be $SL(n)$, for example.
Then we get non-empty equivalued strata from pairs $(w,a/b)$ (with $a/b$
in least common terms and $b=o(w)$) for which $w$ is a power of either an
$n$-cycle or  an $(n-1)$-cycle. 

\section{Admissible subsets of $X(\mathcal O)$}\label{sec.adm.sub} 
Before we describe the structure of the sets $\mathfrak t_w(\mathcal O)_r$
we need a few preliminary remarks and definitions. Consider a scheme $X$ of
finite type over
$\mathcal O$. Then for any positive integer
$N$ Greenberg's functor \cite{greenberg61} (see appendix \ref{app.gre} for
a review)  provides us with a scheme of finite type over $k$ whose set of
$k$-points is
$X(\mathcal O/\epsilon^N
\mathcal O)$. In general the natural $k$-morphism $X(\mathcal
O/\epsilon^{N+1}\mathcal O) \to X(\mathcal O/\epsilon^N \mathcal O)$ can be
complicated, but when $X$ is smooth over $\mathcal O$, as we will always
assume in this section, Greenberg
\cite{greenberg63} shows that 
$X(\mathcal
O/\epsilon^{N+1}\mathcal O)$ is an affine space bundle over  $X(\mathcal
O/\epsilon^N
\mathcal O)$ (more precisely, a torsor under the pullback to $X(\mathcal
O/\epsilon^N \mathcal O)$ of the tangent bundle on $X(k)$). In particular
each $X(\mathcal O/\epsilon^N \mathcal O)$ is then smooth over $k$, and  the
mapping $X(\mathcal
O/\epsilon^{N+1}\mathcal O) \to X(\mathcal O/\epsilon^N \mathcal O)$ is
open and surjective. 

In case $X$ is affine $n$-space $\mathbb A^n$ over $\mathcal O$ (for
example, $\mathfrak t_w(\mathcal O)$ or
$\mathbb A(\mathcal O)$, the two main cases of interest in this
paper), the situation is particularly simple, since then
$X(\mathcal O/\epsilon^N
\mathcal O)$ is $\mathbb A_k^{nN}$, and $X(\mathcal
O/\epsilon^{N+1}\mathcal O) \to X(\mathcal O/\epsilon^N \mathcal O)$ is a
projection map $\mathbb A_k^{n(N+1)} \to \mathbb A_k^{nN}$. 

\subsection{Admissible subsets}
For any positive integer $N$ we write $p_N:X(\mathcal O) \to X(\mathcal
O/\epsilon^N \mathcal O)$ for the canonical surjection (induced of course
by $\mathcal O \twoheadrightarrow \mathcal
O/\epsilon^N \mathcal O$). We say that a subset $Y$ of $X(\mathcal O)$ is
\emph{$N$-admissible} if $Y=p_N^{-1}p_NY$, in which case we introduce $Y_N$
as a convenient notation for $p_NY$. If $Y$ is $N$-admissible, it is clear
that $Y$ is $N'$-admissible for all $N' \ge N$. We say that $Y$ is
\emph{admissible} if there exists $N$ such that $Y$ is $N$-admissible.

\subsection{Topological notions for admissible subsets}
Let $Y$ be an admissible subset of $X(\mathcal O)$. We say that $Y$ is
\emph{open} (respectively, \emph{closed, locally closed, irreducible}) in
$X(\mathcal O)$ if
$Y_N$ is open (respectively, closed, locally closed, irreducible) in
$X(\mathcal O/\epsilon^N \mathcal O)$ for some (equivalently, every)
positive integer $N$ such that $Y$ is $N$-admissible.  (To see the
equivalence of ``some'' and ``every''  use Lemma
\ref{lem.top.adm}.) We define the \emph{closure} $\bar Y$ of $Y$ in
$X(\mathcal O)$ as follows: choose $N$ such that $Y$ is $N$-admissible and
put $\bar Y=p_N^{-1}\bar Y_N$, where $\bar Y_N$ of course denotes the
closure of $Y_N$ in $X(\mathcal O/\epsilon^N \mathcal O)$; by Lemma
\ref{lem.adm.cl}(3) $\bar Y$ is independent of the choice of $N$.

\subsection{Non-singularity for admissible subsets}
Now assume that $Y$ is a locally closed admissible subset of $X(\mathcal
O)$. For $N$ such that $Y$ is $N$-admissible we regard $Y_N$ as a reduced
scheme of finite type over $k$ by putting the induced reduced subscheme
structure on the locally closed subset $Y_N \subset X(\mathcal O/\epsilon^N
\mathcal O)$. We say that $Y$ is \emph{non-singular} if $Y_N$ is
non-singular for some (equivalently, every) positive integer $N$ such that
$Y$ is $N$-admissible. (To see the
equivalence of ``some'' and ``every''  use Lemma
\ref{lem.adm.sm}.) 

Now assume that $Y$ is indeed non-singular, and let $y \in Y$. We want to
define the tangent space $T_{Y,y}$ to $Y$ at $y$. This will be an
admissible $k$-linear subspace of the tangent space $T_{X(\mathcal O),y}$
(see \ref{sub.gre.tan}) to $X(\mathcal O)$ at $y$. 

To define $T_{Y,y}$ we choose $N$ so large that $Y$ is $N$-admissible, and
we denote by $\bar y$ the image of $y$ in $Y_N$. Then the tangent space 
$T_{Y_N,\bar y}$ is a linear subspace of the tangent space $T_{X(\mathcal
O/\epsilon^N\mathcal O),\bar y}$. Recall from \ref{sub.gre.tan} that there
is a canonical surjection 
\begin{equation}\label{eq.sur.tan} 
T_{X(\mathcal O),y} \twoheadrightarrow T_{X(\mathcal
O/\epsilon^N\mathcal O),\bar y},
\end{equation} 
which identifies $T_{X(\mathcal
O/\epsilon^N\mathcal O),\bar y}$ with $T_{X(\mathcal O),y} \otimes_\mathcal
O (\mathcal
O/\epsilon^N\mathcal O)$. We now define $T_{Y,y}$ to be the inverse image
under \eqref{eq.sur.tan} of $T_{Y_N,\bar y}$. It is easy to see that
$T_{Y,y}$ is  independent of the choice of $N$. 

\subsection{Smoothness for maps between admissible
subsets}\label{sub.sm.ad.s} 
 
Now let $f:X \to X'$ be an $\mathcal O$-morphism between smooth
schemes $X,X'$ over $\mathcal O$. For each positive integer $N$
Greenberg's functor yields a $k$-morphism 
\[
f_N:X(\mathcal O/\epsilon^N \mathcal O) \to X'(\mathcal
O/\epsilon^N \mathcal O).
\] 

Suppose that $Y,Y'$ are admissible locally closed subsets of
$X(\mathcal O)$, $X'(\mathcal O)$ respectively, with the
property that $f(Y) \subset Y'$, and let $g:Y \to Y'$ denote the
map obtained by restriction from $f:X(\mathcal O) \to
X'(\mathcal O)$. For each $N$ such that both $Y,Y'$ are
$N$-admissible, we obtain (by restriction from $f_N$) a
$k$-morphism 
\[
g_N:Y_N \to Y'_N.
\] 
As usual we put the induced reduced subscheme structures on
$Y_N,Y_N'$. 

For $M\ge N$ there is a commutative square 
\[
\begin{CD}
Y_M @>g_M>> Y'_M \\
@V VV @VVV \\
Y_N @>g_N>> Y'_N
\end{CD}
\]
in which the vertical arrows (the obvious surjections) are smooth
 (by Lemma
\ref{lem.adm.sm}). It then follows from EGA IV (17.11.1) that 
if $g_M$ is smooth, then $g_N$ is also smooth. However, if $g_N$
is smooth, it is not necessarily the case that $g_M$ is smooth.

We say that $Y$ is \emph{smooth} over $Y'$ (or that $g:Y \to Y'$
is \emph{smooth}) if
$g_M$ is smooth for all $M \ge N$. The remarks we just made show
that this condition is independent of the choice of $N$ for which
$Y,Y'$ are both $N$-admissible. 

It is evident from the definitions that if $g:Y \to Y'$ is
smooth, and $Y'$ is non-singular, then $Y$ is non-singular.
Using  Lemma \ref{lem.adm.sm} one checks easily that if $Y$ is
smooth over $Y'$, then $g^{-1}Z'$ is smooth over $Z'$ for any
admissible locally closed subset $Z'$ of $X'(\mathcal O)$ such
that $Z' \subset Y'$. 

\begin{lemma}\label{lem.sm.dif} Now assume that 
both $Y,Y'$ are  non-singular. Define the differential
$dg_y$ of $g:Y \to Y'$ at $y \in Y$ to be the $k$-linear map
$dg_y:T_{Y,y} \to T_{Y',f(y)}$    obtained by restricting  the
differential
$df_y:T_{X(\mathcal O),y}
\to T_{X'(\mathcal O),f(y)}$ to the tangent space $T_{Y,y}$.  
Then $Y$ is smooth over $Y'$ if
and only if the differential 
\begin{equation}\label{eq.sur.Y}
dg_y:T_{Y,y} \to T_{Y',f(y)}
\end{equation}
is surjective for all $y \in Y$. 
\end{lemma}
\begin{proof}
Let $y \in Y$ and put $y':=f(y)$.  For any integer $M$ with $M
\ge N$ we denote by $y_M$ the image of $y$ under the canonical
surjection $Y \twoheadrightarrow Y_M$. We do the same for $y'$,
so that $y_M'=g_M(y_M)$. We then have a commutative square 
\begin{equation}\label{sq.tan.M}
\begin{CD}
T_{Y,y} @>>> T_{Y',y'} \\
@V VV @VVV \\
T_{Y_M,y_M} @>>> T_{Y'_M,y'_M}
\end{CD}
\end{equation}
 in which the horizontal maps are differentials and the vertical
maps are the canonical surjections. 

The implication ($\Longleftarrow$) of the lemma is now clear,
since \eqref{sq.tan.M} together with the surjectivity of
\eqref{eq.sur.Y} shows that each $g_M$ is a submersion. 

It remains to prove the reverse implication ($\Longrightarrow$).
To simplify notation we put $L:=T_{X(\mathcal O),y}$ (a free
$\mathcal O$-module of finite rank) and $V:=T_{Y,y}$ (an
$N$-admissible $k$-linear subspace of $L$), and we use parallel
notation for $Y'$. Then  \eqref{sq.tan.M} becomes the square 
\[
\begin{CD}
V @>\psi>> V' \\
@V VV @VVV \\
V/\epsilon^M L @>\psi_M>> V'/\epsilon^M L',
\end{CD}
\]
where $\psi$ is the $k$-linear map obtained by restriction from
the $\mathcal O$-linear map $\varphi:L \to L'$ defined by
$\varphi:=df_y$.

Our assumption that $Y$ is smooth over $Y'$ tells us that
$\psi_M $ is surjective for all $M \ge N$, which just means that 
\begin{equation}\label{eq.phiV}
\varphi V+\epsilon^M L' =V'
\end{equation}
for all $M \ge N$. Since $V' \supset \epsilon^NL'$, we conclude
that 
\[
\epsilon^NL' \subset \bigcap_{M \ge N}(\varphi L +\epsilon^M L')
=\varphi L. 
\]
(Here we used that $\varphi L$ is an $\mathcal O$-submodule of
$L'$.) Since $V \supset \epsilon^N L$, we see that 
\[
\varphi V \supset \epsilon^N  \varphi L \supset \epsilon^{2N}L'.
\] 
Taking $M=2N$ in \eqref{eq.phiV}, we conclude that $\varphi
V=V'$, showing that \eqref{eq.sur.Y} is surjective, as desired. 
\end{proof} 
\begin{corollary}
Suppose that $Y$ is smooth over $Y'$. Then for each $y \in Y$
the $F$-linear map obtained by extension of scalars from the
differential 
\begin{equation}\label{eq.nd.label}
df_y:T_{X(\mathcal O),y} \to T_{X'(\mathcal O),f(y)}
\end{equation}
is surjective.  In other words, $Y$ is necessarily contained in
the subset of $X(F)$ consisting of all points at which the 
$F$-morphism obtained by extension of scalars from $f:X \to X'$
is smooth. 
\end{corollary}
\begin{proof}
Put $y':=f(y)$. Choose $N$ for which both $Y,Y'$ are
$N$-admissible. Write $y'_N$ for the image of $y'$ under the
canonical surjection $Y' \twoheadrightarrow Y'_N$, and let $Z'$
denote the preimage of $y'_N$ in $Y'$.  We have noted before
that $g^{-1}Z'$ is smooth over $Z'$. Since $Z'$ is obviously
non-singular, so too is $g^{-1}Z'$. The previous lemma, applied
to $g^{-1}Z' \to Z'$,  then tells us that the map 
\begin{equation*}
T_{g^{-1}Z',y} \to T_{Z',y'}
\end{equation*}
is surjective, and hence that the image of the map
\eqref{eq.nd.label} contains $\epsilon^N T_{X'(\mathcal
O),y'}$.    Therefore
the $F$-linear map obtained from  
\eqref{eq.nd.label} is surjective. 
\end{proof}

\subsection{Codimensions of admissible subsets}
Now suppose  that $X(k)$ is irreducible (which implies that $X(\mathcal
O/\epsilon^N \mathcal O)$ is irreducible for every positive integer~$N$).  
When
$Y$ is an admissible, locally closed,  irreducible subset of $X(\mathcal
O)$, we define its
\emph{codimension} in $X(\mathcal O)$ to be the codimension of $Y_N$ in
$X(\mathcal O/\epsilon^N \mathcal O)$ for any $N$ such that $Y$ is
$N$-admissible; it is easy to see that this notion of codimension is
independent of the choice of $N$.

\section{Structure of the strata $\mathfrak t_w(\mathcal  O)_r$} 
We now continue the discussion of $\mathfrak t_w(\mathcal O)_r$, retaining 
the notation used before. In particular $l$ denotes the order of $w$, and
$\mathfrak a_m$ denotes   the linear subspace  defined in subsection
\ref{sub.wh.ne}. We now have the right vocabulary to discuss the structure
of
$\mathfrak t_w(\mathcal O)_r$. 

In the next result we will need the stabilizer $W_{w,r}$ of the
pair
$(w,r)$; thus $W_{w,r}$ consists of elements $x \in W$ such that $xw=wx$
and $xr=r$. We will also need the eigenspaces 
\[
(\mathfrak
t/\mathfrak a_{j+1})(w,j):=\{v \in (\mathfrak
t/\mathfrak a_{j+1}): wv=\zeta_l^{-j}v \}.
\]

\begin{proposition}\label{prop.t.str}
Assume that $Y:=\mathfrak t_w(\mathcal O)_r$ is non-empty, which guarantees
in particular that $r$ takes values in $\frac{1}{l}\mathbb Z$. Let $N$ be a
positive integer large enough that $r(\alpha) < N$ for all $\alpha \in R$. 
Then the following conclusions hold.
\begin{enumerate}
\item The subset $Y$ of $\mathfrak t_w(\mathcal O)$ is $N$-admissible. Thus
$Y$ is the preimage of its image $Y_N$ in $\mathfrak t_w(\mathcal
O/\epsilon^N \mathcal O)$. 
\item The closure $\bar Y$ of $Y$ in $\mathfrak t_w(\mathcal O)$ is the
admissible
$k$-linear subspace
\[
\{ u \in \mathfrak t_w(\mathcal O): \val \alpha(u) \ge r(\alpha) \quad
\forall \, \alpha \in R \}
\]
of $\mathfrak t_w(\mathcal O)$. 
\item The subset $Y \subset \bar Y$ is the complement of finitely many
admissible $k$-linear hyperplanes $H_1,\dots,H_m$ in $\bar Y$.
Consequently $Y$ is locally closed, irreducible, and non-singular. The
group $W_{w,r}$ preserves $Y$, $\bar Y$ and permutes the hyperplanes 
$H_1,\dots,H_m$ in $\bar Y$; moreover, $W_{w,r}$ acts freely on $Y_N$. 
\item The codimension of $Y$ in $\mathfrak t_w(\mathcal O)$ is the same as
that of the linear subspace $\bar Y$, namely 
\[
\dim_k(\mathfrak t_w(\mathcal O)/\bar Y)=\sum_{j=0}^\infty \dim_k (\mathfrak
t/\mathfrak a_{j+1})(w,j).
\]
When $w=1$, this expression for the codimension simplifies to
\[\sum_{j=0}^\infty j \cdot\dim_k (\mathfrak a_{j+1}/\mathfrak a_j).
\] 
\end{enumerate}
\end{proposition}
\begin{proof}
(1) Suppose that $u \in \mathfrak t_w(\mathcal O)_r$ and that $u' \in
\epsilon^N \mathfrak t_w(\mathcal O)$. For each root $\alpha$ we must show
that $\alpha(u+u')$ has the same valuation as $\alpha(u)$. This is clear
from our hypothesis that  $r(\alpha) < N$. 

(2) By Lemma \ref{lem.t.str.ne} the non-emptiness of $\mathfrak
t_w(\mathcal O)_r$ implies the non-emptiness of $\mathfrak
t(w,j)\cap \mathfrak a_{j+1}^\sharp$ for all $j\ge 0$. Since 
$\mathfrak
t(w,j)\cap \mathfrak a_{j+1}^\sharp$ is the complement of finitely many
hyperplanes in $\mathfrak
t(w,j)\cap \mathfrak a_{j+1}$, we see that the closure of $\mathfrak
t(w,j)\cap \mathfrak a_{j+1}^\sharp$ is $\mathfrak
t(w,j)\cap \mathfrak a_{j+1}$. Lemma \ref{lem.t.str.ne} then implies that 
the closure of
$\mathfrak t_w(\mathcal O)_r$ is the set of all $u=\sum_{j=0}^\infty u_j 
\epsilon_E^j
\in \mathfrak t_w(\mathcal O)$ such that $u_j \in  \mathfrak a_{j+1}$
for all $j \ge 0$.   On the other hand $u$ lies in 
\[
\{ u \in \mathfrak t_w(\mathcal O): \val \alpha(u) \ge r(\alpha) \quad
\forall \, \alpha \in R \}
\]
if and only if $\alpha(u_j)=0$ whenever $j/l < r(\alpha)$, and this happens
if and only if $u_j \in \mathfrak a_{j+1}$ (by the very definition of
$\mathfrak a_{j+1}$). 

(3) The proof of (2) shows that $Y$ is the complement of finitely many
admissible hyperplanes in $\bar Y$.  
The freeness of the action of $W_{w,r}$ on $Y_N$ follows from the freeness 
(see \ref{sub.free}) 
of the action of $W_w$ on the larger set $\mathfrak t_w(\mathcal
O/\epsilon^N\mathcal O)_{r<N}$. The remaining
statements are clear. 

(4) The description of $\bar Y$ given in (2) shows that 
\[
\dim_k(\mathfrak t_w(\mathcal O)/\bar Y)=\sum_{j=0}^\infty \dim_k \mathfrak
t(w,j)/(\mathfrak
t(w,j) \cap \mathfrak a_{j+1}).
\]
Since the order of $w$ is invertible in $k$, we see that $\mathfrak
t(w,j)/(\mathfrak
t(w,j) \cap \mathfrak a_{j+1})$ can be identified with $(\mathfrak
t/\mathfrak a_{j+1})(w,j)$. 

Finally, when $w=1$, we have $(\mathfrak
t/\mathfrak a_{j+1})(w,j)=\mathfrak
t/\mathfrak a_{j+1}$, whose dimension is $\sum_{j' =
j+1}^\infty\dim_k(\mathfrak a_{j'+1}/\mathfrak a_{j'})$. This proves the
last statement in (4). 
\end{proof}

\section{Strata in $\mathbb A(\mathcal O)'$} 
We now stratify $\mathbb A(\mathcal O)'$. We obtain the desired strata in 
$\mathbb A(\mathcal O)'$ as images of the strata $\mathfrak t_w(\mathcal
O)_r$ that we have already studied. 

\subsection{Definition of the map $f_w:\mathfrak t_w(\mathcal O) \to \mathbb
A(\mathcal O)$} \label{sub.def.fw} 
Let $w \in W$. 
The map $\mathfrak t(\mathcal O_E) \to \mathbb A(\mathcal O_E)$ on
$\mathcal O_E$-points induced by our morphism $ f$ restricts to a
map 
\[
f_w:\mathfrak t_w(\mathcal O) \to \mathbb A(\mathcal O).
\]
Recall that the centralizer $W_w$ (of $w$ in $W$) acts on $\mathfrak
t_w(\mathcal O)$. This action preserves the fibers of the map $f_w$. 

In fact $f_w$ comes from a morphism of schemes over
$\mathcal O$ that will also be denoted simply by $f_w:\mathfrak t_w \to
\mathbb A_{\mathcal O}$, with $\mathbb A_{\mathcal O}$ denoting the
$\mathcal O$-scheme obtained from $\mathbb A$ by extending scalars from
$k$ to $\mathcal O$. This is best understood using the point of view (see
\ref{sub.tw.fp}) that
$\mathfrak t_w$ is the fixed point scheme of a $\mathbb Z/l\mathbb Z$-action
on
$R_{\mathcal O_E/\mathcal O}\mathfrak t$. (We again remind the reader that
fixed point schemes are discussed in appendix \ref{sec.fp.gx}.)

By $R_{\mathcal O_E/\mathcal O}\mathbb A$ we will of course mean the scheme
obtained from $\mathbb A$ by extending scalars from $k$ to $\mathcal O_E$,
and then (Weil) restricting scalars from $\mathcal O_E$ to $\mathcal O$.
Thus 
\[
(R_{\mathcal O_E/\mathcal O}\mathbb A)(A)=\mathbb A(A \otimes _\mathcal O
\mathcal O_E)
\]
for any $\mathcal O$-algebra $A$.

As in subsection \ref{sub.tw.fp} the automorphism $\tau_E$ of $\mathcal
O_E/\mathcal O$ induces an automorphism $\tau_E$ of $R_{\mathcal
O_E/\mathcal O}\mathbb A$ of order $l$, so that we obtain an action of
$\mathbb Z/l\mathbb Z$ on $R_{\mathcal O_E/\mathcal O}\mathbb A$. The fixed
point scheme of $\mathbb Z/l\mathbb Z$  on $R_{\mathcal O_E/\mathcal
O}\mathbb A$ is $\mathbb A_{\mathcal O}$, as one sees from the (easy) fact
that
$A$ is the set of fixed points of $\id_A \otimes \tau_E$ on $A
\otimes_{\mathcal O} \mathcal O_E$ for any $\mathcal O$-algebra (or even
$\mathcal O$-module) 
$A$. 

Starting from $f:\mathfrak t \to \mathbb A$, then extending scalars to
$\mathcal O_E$, then restricting scalars to $\mathcal O$, we get an
$\mathcal O$-morphism 
\[
R(f):R_{\mathcal O_E/\mathcal O}\mathfrak t \to R_{\mathcal O_E/\mathcal
O}\mathbb A 
\]
which intertwines the automorphism $w\tau_E$ of $R_{\mathcal O_E/\mathcal
O}\mathfrak t$ with the automorphism $\tau_E$ of $R_{\mathcal O_E/\mathcal
O}\mathbb A$, and hence induces the desired $\mathcal O$-morphism $f_w$
upon taking fixed points under $\mathbb Z/l\mathbb Z$. 

\subsection{Fibers of $\mathfrak t_w(\mathcal O)' \to \mathbb A(\mathcal
O)'$} \label{sub.fib.TWA} 
Recall that the centralizer $W_w$ acts on $\mathfrak t_w(\mathcal O)$,
preserving the fibers of $f_w$. We claim that $W_w$ acts simply
transitively on every non-empty fiber of the restriction of $f_w$ to
$\mathfrak t_w(\mathcal O)'$. Indeed, let $u,u' \in \mathfrak t_w(\mathcal
O)'$ and suppose that $f_w(u)=f_w(u')$. Then there exists unique $x \in W$
such that $xu=u'$. Using that $u$, $u'$ are fixed by $w\tau_E$, we see that 
\[
wx\tau_E(u)=u'=xw\tau_E(u)
\]
and hence that $wx=xw$, as claimed. 

\subsection{Definition of the strata $\mathbb A(\mathcal O)_s$ in 
$\mathbb A(\mathcal O)'$} 
Consider a pair $(w,r) \in 
W \times \mathcal R$. We denote by $\mathbb A(\mathcal O)_{w,r}$ the image
of $\mathfrak t_w(\mathcal O)_r$ under the map $f_w$. It is clear that 
$\mathbb A(\mathcal O)_{w,r}$ depends only on the $W$-orbit of $(w,r)$
(with, as usual, $W$ acting on itself by conjugation). Thus it is
often better  to index the strata by the set $\mathcal S$ of orbits of $W$
on
$W \times \mathcal R$. In other words, given  $s \in \mathcal S$,
represented by a pair $(w,r)$, we will often write 
$\mathbb A(\mathcal O)_s$ instead of $\mathbb A(\mathcal O)_{w,r}$.

Since $\mathbb A(\mathcal O)_s$ is by definition obtained as the image of
$\mathfrak t_w(\mathcal O)_r$, Proposition \ref{prop.ne.gen} tells us when 
$\mathbb A(\mathcal O)_s$ is non-empty. 

\begin{lemma}\label{lem.weyl.elt}
The set $\mathbb A(\mathcal O)'$ is the disjoint union of the strata
$\mathbb A(\mathcal O)_s$. 
\end{lemma}   
\begin{proof}
Let $c \in \mathbb A(\mathcal O)'$. Consider the fiber over $c$ 
 of the map $\mathfrak t(\bar F) \to \mathbb A(\bar
F)$ induced by our morphism $f:\mathfrak t \to \mathbb A$. (Recall that for
any
$k$-algebra $A$ we have  
$\mathfrak t(A)=\mathfrak t \otimes_k A$.) The Weyl group $W$ acts simply
transitively on this fiber. Moreover, 
 since $f_{\reg}$ is \'etale, the fiber is actually contained
in the subset $\mathfrak t(F_{\sep})$ of $\mathfrak t(\bar F)$. 
The action of $\Gal(F_{\sep}/F)$ on $\mathfrak t(F_{\sep})$ preserves the
fiber because $c$ is defined over $F$ (and even over $\mathcal O$). 

Now choose
an element $u$ in the fiber.  For any element $\tau \in \Gal(F_{\sep}/F)$
there exists a unique
$w_\tau
\in W$ such that $w_\tau \tau(u)=u$, and $\tau \mapsto w_\tau$ is a
homomorphism from $\Gal(F_{\sep}/F)$ to $W$. Since $|W|$ is invertible in
$k$, this homomorphism factors through the quotient $\Gal(F_{\tame}/F)$ of
$\Gal(F_{\sep}/F)$, and in fact we will now simply regard $\tau\mapsto
w_\tau$ as a homomorphism from $\Gal(F_{\tame}/F)$ to $W$. Recall from
before the topological generator $\tau_\infty$ of $\Gal(F_{\tame}/F)$.
Putting $w:=w_{\tau_\infty}$, we have associated an element $w \in W$ to
the element $u$ in the fiber. As usual we write $l$ for $o(w)$ and $E$ for
$F_l$.

It is clear from the definitions that $u \in \mathfrak t_w( F)$ and
that $u \mapsto c$ under our morphism $f$.  
The valuative criterion of properness,
applied to the proper morphism $f$ and the valuation ring $\mathcal O_E$,
implies that $u \in \mathfrak t(\mathcal O_E)$ and hence that $u \in
\mathfrak t_w(\mathcal O)$. 
 Define $r \in \mathcal R$ by
$r(\alpha):=\val
\alpha(u)$. Then $u$ lies in the stratum $\mathfrak t_w(\mathcal O)_r$,
and therefore $c$ lies in the stratum
$\mathbb A(\mathcal O)_s$. Thus we have shown that our strata exhaust
$\mathbb A(\mathcal O)$. 

It remains to establish disjointness of our strata. 
Suppose that $u_1 \in \mathfrak t_{w_1}(\mathcal O)_{r_1}$ and that 
$u_2 \in \mathfrak t_{w_2}(\mathcal O)_{r_2}$. Suppose further that
$u_1$ and $u_2$ have the same image $c$ in $\mathbb A(\mathcal O)$.  We must
show that $(w_1,r_1)$ and $(w_2,r_2)$ are in the same $W$-orbit. This is
easy: there exists a unique element $x \in W$ such that $xu_1=u_2$, and
this element $x$ transforms $(w_1,r_1)$ into $(w_2,r_2)$. 
\end{proof}

\begin{lemma}\label{lem.inv.i.as}
We have 
\[
f_w^{-1}(\mathbb A(\mathcal O)_s)=\coprod_{x \in W_w/W_{w,r}}\mathfrak
t_w(\mathcal O)_{xr}, 
\]
where, as usual, $W_{w,r}$ denotes the stabilizer in $W$ of the pair
$(w,r)$, and $s$ denotes the $W$-orbit of $(w,r)$. Moreover, $W_{w,r}$ acts
simply transitively on each fiber of $\mathfrak t_w(\mathcal O)_r
\twoheadrightarrow \mathbb A(\mathcal O)_s$. 
\end{lemma}
\begin{proof}
This follows from the discussion in subsection \ref{sub.fib.TWA} and the
obvious equality (valid for any $x \in W_w$)
\[
\mathfrak t_w(\mathcal O)_{xr}=x\mathfrak t_w(\mathcal O)_{r}.
\]
\end{proof}

\section{Structure of the strata $\mathbb A(\mathcal O)_s$ in $\mathbb
A(\mathcal O)'$} 
In this section, after introducing a couple of definitions, we are going to
formulate Theorem \ref{thm.a.str}, which describes the structure of the
strata $\mathbb A(\mathcal O)_s$ in $\mathbb
A(\mathcal O)'$. 
\subsection{Definitions of $\delta_r$ and $c_w$} 
Let $s \in S$ be the $W$-orbit of  the pair $(w,r)$. Let us assume that the
stratum $\mathbb A(\mathcal O)_s$ is non-empty (equivalently: 
$\mathfrak t_w(\mathcal O)_r$ is non-empty). As usual we put $l:=o(w)$ and
$E:=F_l$. 

So far we have not used the $F$-torus $T_w$ that goes along with $\mathfrak
t_w(F)$.  This torus splits over $E$ and is obtained by using $w$ to twist 
the (split) torus over $F$ obtained by extension of scalars from $T$. In
particular we have 
\[
T_w(F)=\{t \in T(E):w\tau_E(t)=t \},
\] 
and the Lie algebra of $T_w$ is canonically isomorphic to $\mathfrak
t_w(F)$. 

There is a canonical $G(F)$-conjugacy class of $F$-embeddings  $T_w \to
G$ (with image a maximal $F$-torus in $G$). This is well-known (perhaps see
\cite{gkm.pre2} for a rather concrete presentation of this material). Fixing
such an embedding, we may identify $\mathfrak t_w(F)$ with a Cartan
subalgebra in $\mathfrak g(F)$. For any regular element $u \in \mathfrak
t_w(F)$ the centralizer in $\mathfrak g(F)$  of $u$  is equal to $\mathfrak
t_w(F)$, and we have the usual non-zero scalar $\Delta(u)$ in $F$ defined
by 
\[
\Delta(u):=\det(\ad(u);\mathfrak g(F)/\mathfrak t_w(F)).
\]
Clearly this determinant is simply the product of the values on $u$ of all
the roots of our Cartan subalgebra. Therefore, if $u \in \mathfrak
t_w(\mathcal O)_r$, we have 
\[
\val \Delta(u) =\delta_r,
\]
where 
\[
\delta_r:=\sum_{\alpha \in R}r(\alpha).
\]
Note that $\delta_r$ is a non-negative integer. (It is clearly
non-negative, and our expression for it as the valuation of 
$\Delta(u) \in F^\times$ shows that it is an integer. We could have defined
$\delta_r$  without this digression concerning $T_w$,
but then it would not have been clear that $\delta_r$ is an integer.)

Since $r(-\alpha)=r(\alpha)$ for all $\alpha \in R$ (because of our
assumption that $\mathfrak t_w(\mathcal O)_r$ be non-empty), 
we also have 
\[
\delta_r=2\sum_{\alpha \in R^+}r(\alpha).
\]

We need some more notation before we state the next result. We denote by
$\mathfrak t^w$ the set of fixed points of $w$ on $\mathfrak t$, and we
denote by $c_w$ the integer 
\[
c_w:=\dim_k\mathfrak t -\dim_k \mathfrak t^w.
\]
Equivalently (because of our hypothesis on the characteristic of our base
field), $c_w$ is the dimension of $T$ minus the dimension of the maximal
$F$-split torus in $T$. 

\subsection{Valuation of the Jacobian of $f_w$} 
Recall  the map 
$f_w:\mathfrak t_w(\mathcal O) \to \mathbb A(\mathcal O)$, which, as we
saw in  \ref{sub.def.fw}, comes from a morphism of schemes over
$\mathcal O$. Our chosen basic invariants allow us to identify $\mathbb A$
with $\mathbb A^n$, and by choosing an $\mathcal O$-basis of the free
$\mathcal O$-module $\mathfrak t_w(\mathcal O)$, we may also identify the
$\mathcal O$-scheme $\mathfrak t_w$ with
$\mathbb A^n_\mathcal O$.

These identifications allow us to think of the differential $(df_w)_u$ of
$f_w$ at $u \in \mathfrak t_w(\mathcal O)$ concretely as a square matrix
$D_u
\in M_n\mathcal O$, as in \ref{sub.a.p.d}. In the next lemma we will compute
the valuation of $\det D_u$ for $u \in
\mathfrak t_w(\mathcal O)_r$, as this will be needed in the proof of
Theorem \ref{thm.a.str}. Observe that making a different choice of
$\mathcal O$-basis for
$\mathfrak t_w(\mathcal O)$ does not affect the valuation of $\det D_u$, so
that it makes sense to write $\val\det (df_w)_u$.

\begin{lemma}\label{lem.val.Jac} 
For any $u \in \mathfrak t_w(\mathcal O)_r$ the non-negative integer $\val
\det (df_w)_u$ is equal to $(\delta_r+c_w)/2$.
\end{lemma}
\begin{proof}
We can
calculate this determinant after extending scalars from $\mathcal O$  to
$\mathcal O_E$. Then we are dealing with the $\mathcal O_E$-linear map 
\[
\id_E \otimes (df_w)_u : \mathcal O_E \otimes_\mathcal O \mathfrak
t_w(\mathcal O) \to \mathcal O_E^n,
\]
which is none other than the restriction of 
\[(df)_u:\mathfrak t(\mathcal O_E) \to \mathcal O_E^n
\]
to the subspace $\mathcal O_E \otimes_\mathcal O \mathfrak t_w(\mathcal O)$
of $\mathfrak t(\mathcal O_E)$. We conclude that 
\[
\val\det (df_w)_u=\val\det(df)_u + \frac{1}{l}\dim_k \frac{\mathfrak
t(\mathcal O_E)}{\mathcal O_E \otimes_\mathcal O \mathfrak t_w(\mathcal O)}.
\]

In order to prove the lemma it is enough to check that 
\begin{equation}\label{eq.help1}
\val\det(df)_u=\delta_r/2 
\end{equation} 
and that 
\begin{equation}\label{eq.help2}
\dim_k \frac{\mathfrak t(\mathcal
O_E)}{\mathcal O_E \otimes_\mathcal O \mathfrak t_w(\mathcal O)}=lc_w/2.
\end{equation}
Now \eqref{eq.help1} follows from \eqref{eq.jac.pro}, and \eqref{eq.help2}
is \cite[Lemma 3]{bezrukavnikov96}.
(Bezrukavnikov treats simply connected groups over $\mathbb C$, but his
proof goes through in our situation. For this we just need to show that the
representation of $W$ on $\mathfrak t$ is isomorphic to its own
contragredient. Since $|W|$ is invertible in $k$, 
it is enough to check that the $W$-modules $\mathfrak t$ and
$\mathfrak t^*$ have the same character, and this is clear, since
$\mathfrak t$ is obtained by tensoring
$X_*(T)$ with $k$, so that all character values lie in the prime field.)
\end{proof}

Let us introduce one more bit of notation before stating the next theorem.
We put $e(w,r):=(\delta_r+c_w)/2$. It follows from Lemma \ref{lem.val.Jac}
that $e(w,r)$ is an integer (non-negative, of course). 

\begin{theorem}\label{thm.a.str} Let $s \in S$ be the $W$-orbit of  the pair $(w,r)$. Let us assume that the
stratum $\mathbb A(\mathcal O)_s$ is non-empty. Then we have the following
conclusions.
\begin{enumerate}
\item The subset $\mathbb A(\mathcal O)_s$ of $\mathbb A(\mathcal O)$ is
admissible; more precisely, it is $N$-admissible whenever $N>2e(w,r)$. 
Moreover it is locally closed, irreducible and non-singular. 
\item The codimension of $\mathbb A(\mathcal O)_s$ in $\mathbb A(\mathcal
O)$ is given by 
\begin{equation}
d(w,r) + e(w,r),
\end{equation}
where $d(w,r)$ is the codimension of $\mathfrak t_w(\mathcal O)_r$ in
$\mathfrak t_w(\mathcal O)$. 
\item $\mathfrak t_w(\mathcal O)_r$ is smooth over $\mathbb A(\mathcal
O)_s$. Here we are using the notion of smoothness discussed in
\textup{\ref{sub.sm.ad.s}}. 
\end{enumerate}
\end{theorem}
\begin{proof}
This will be proved in  section \ref{sec.pr.m.th}. 
\end{proof}

\section{Relation between the strata $\mathfrak t_w(\mathcal O)_r$ and
$\mathbb A(\mathcal O)_s$}
\subsection{Set-up for this section} \label{sub.setup91}
Consider a non-empty stratum $\mathfrak
t_w(\mathcal O)_r$, and let $s$ again denote the $W$-orbit of the pair
$(w,r)$. We abbreviate $e(w,r)$ to $e$. We now have a good
understanding of the strata
$\mathfrak t_w(\mathcal O)_r$ and
$\mathbb A(\mathcal O)_s$, but we would  like to supplement this by
analyzing the smooth morphism  
\[
\mathfrak t_w(\mathcal O)_r \to \mathbb A(\mathcal O)_s 
\] 
obtained by restriction from $f_w$. 

To do so we choose $N$ large enough that $N > 2e$, and use the fact that
both $\mathfrak t_w(\mathcal O)_r$ and  $\mathbb A(\mathcal O)_s$ are
$N$-admissible (see Proposition \ref{prop.t.str} and Theorem
\ref{thm.a.str}). (For additional details see the first few lines of the
proof of Theorem
\ref{thm.a.str}, where it is shown that $\mathfrak t_w(\mathcal O)_r$  is
even $(N-e)$-admissible.) In Theorem \ref{thm.HH} we will gain an
understanding of  
\[
\mathfrak t_w(\mathcal O)_r \to \mathbb A(\mathcal O)_s 
\] 
by analyzing the smooth morphism 
\[
\mathfrak t_w(\mathcal O/\epsilon^N\mathcal O)_r  \to \mathbb A(\mathcal
O/\epsilon^N\mathcal O)_s .
\]
Here we have written $\mathfrak t_w(\mathcal O/\epsilon^N\mathcal O)_r$
for the image of $\mathfrak t_w(\mathcal O)_r$ in $\mathfrak
t_w(\mathcal O/\epsilon^N\mathcal O)$. Similarly we have written  $\mathbb
A(\mathcal O/\epsilon^N\mathcal O)_s$ for the image of $\mathbb
A(\mathcal O)_s$ in $\mathbb A(\mathcal
O/\epsilon^N\mathcal O)$. 

Theorem \ref{thm.HH} makes use of a rank $e$ vector bundle $\tilde V$ 
over $\mathbb A(\mathcal O/\epsilon^N\mathcal O)_s$  that will be
constructed in the course of proving  the theorem. This vector bundle acts
on $\mathfrak t_w(\mathcal O/\epsilon^N\mathcal O)_r$ over $\mathbb
A(\mathcal O/\epsilon^N\mathcal O)_s$. The group $W_{w,r}$ also
acts on $\mathfrak t_w(\mathcal O/\epsilon^N\mathcal O)_r$ over $\mathbb
A(\mathcal O/\epsilon^N\mathcal O)_s$. The two actions commute, so the
group scheme $H:=W_{w,r} \times \tilde V$ (product over $k$) over 
 $\mathbb A(\mathcal O/\epsilon^N\mathcal O)_s$ also acts on 
$\mathfrak t_w(\mathcal O/\epsilon^N\mathcal O)_r$ over $\mathbb
A(\mathcal O/\epsilon^N\mathcal O)_s$.

\begin{theorem} \label{thm.HH}
The space $\mathfrak t_w(\mathcal
O/\epsilon^N\mathcal O)_r$ over $\mathbb
A(\mathcal O/\epsilon^N\mathcal O)_s$ is a torsor under
$H$. In particular we can factorize the morphism $\mathfrak t_w(\mathcal O/
\epsilon^N\mathcal O)_r \to \mathbb A(\mathcal
O/\epsilon^N\mathcal O)_s$ as the composition of two morphisms, one of
which is a bundle of affine spaces of dimension $e$, and the other of
which is an
\'etale covering that is Galois with group $W_{w,r}$. The factorization can
be done in either order. 
\end{theorem}
\begin{proof}
This will be proved in section \ref{sec.pf.H}.
\end{proof}

\section{Behavior of admissibility under polynomial maps 
$f:\mathcal O^n
\to
\mathcal O^n$}

In this section we will establish some technical results needed for the
proofs of the theorems we have stated. It is the morphism $f_w:\mathfrak
t_w(\mathcal O) \to \mathbb A(\mathcal O)$ that we need to understand, but
it is conceptually simpler to work in the more general
context of polynomial maps $f:\mathcal O^n \to \mathcal O^n$. The key Lemma
\ref{lem.suc.app}, a generalization of Hensel's lemma, is a variant of a
special case of one of the main results of  Greenberg's paper
\cite{greenberg66}.

\subsection{Set-up for this section} \label{sub.a.p.d} 
Consider the polynomial ring
$A=\mathcal O[X_1,\dots,X_n]$.  Thus $\Spec A$ is affine
$n$-space $\mathbb A^n$ over $\mathcal O$.   In this section we
study a morphism
$\mathbf f:\mathbb A^n
\to \mathbb A^n$ of schemes over $\mathcal O$. Thus $\mathbf f$
is given by an $n$-tuple $\mathbf f=(f_1,\dots,f_n)$ of elements
in $A$.

We write $L$ for $\mathcal O^n$, the set of $\mathcal O$-valued
points of $\mathbb A^n$.  We are mainly interested in the map
\[
f:L \to L
\]
on $\mathcal O$-valued
points induced by our morphism $\mathbf f:\mathbb A^n
\to \mathbb A^n$.

We regard the differential $d\mathbf f$ of $\mathbf f$
concretely as an element of $M_nA$, the ring of square matrices
 of size $n$ with entries in
$A$. Of course the matrix entries are the partial derivatives 
$\partial f_j/\partial X_i$.  We denote by
$D_x$ the value of
$d\mathbf f$ at
$x
\in L$; thus
$D_x \in M_n\mathcal O$ and $D_x$ can also be viewed as an $\mathcal
O$-linear map $D_x:L \to L$.   It is evident that the reduction modulo
$\epsilon^N$ of
$D_x$ depends only on
$x$ modulo
$\epsilon^N$. (Here $N$ is any non-negative integer.)

We put
$g=\det(d\mathbf f)
\in A$; clearly
$g(x)=\det D_x$ for $x \in L$. The reduction
modulo $\epsilon^N$ of $g(x)$ depends only on $x$ modulo
$\epsilon^N$.
For $x \in L$  we write $d(x)$ for the valuation of $g(x)$. Thus
$d(x)$ is a non-negative integer when $g(x)\ne 0$, and
$d(x)=+\infty$ when $g(x)=0$.   

\subsection{The linear case}\label{sub.lin.case} 
The situation is of course especially simple
when our morphism $f$ is linear. In this subsection we suppose that $f:L \to
L$ is given   by multiplication by a matrix $A \in
M_n\mathcal O$ whose determinant is non-zero, and we put $d:=\val\det A$. 
Thus the $\mathcal O$-module $L/AL$ has length $d$, hence is killed by
$\epsilon^d$, which is to say that 
\begin{equation}\label{eq.edLc}
\epsilon^d L \subset AL.
\end{equation}

\begin{lemma}\label{lem.linear}
 Let $Y$
be a subset of $L$ that is admissible, locally closed, irreducible,
non-singular of codimension $a$ in $L$. Then $AY$ is admissible, 
locally closed, irreducible,
non-singular of codimension $a+d$ in $L$.
\end{lemma}
\begin{proof}
Easy.
\end{proof} 

For any non-negative integer $N$ we denote by $A_N$ the reduction of $A$
modulo $\epsilon^N$. We view $A_N$ as an $\mathcal O/\epsilon^N\mathcal
O$-linear map $A_N:L/\epsilon^N L \to L/\epsilon^N L$. In the next lemma we
will see that for  $N \ge d$ the kernel of $A_N$ is always 
$d$-dimensional and is even independent
of $N$, up to canonical isomorphism.  The kind of canonical isomorphism
that will come up is of the following type: for any integers $M,N$ with
$M\le N$, there is a canonical isomorphism  
$ L/\epsilon^{N-M} L \cong \epsilon^M L /\epsilon^N L$, given of course by
multiplication by $\epsilon^M$.  

\begin{lemma}\label{lem.lin.ker}
 Suppose that $N \ge
d$. Then $\ker(A_N)$ is $d$-dimensional and is contained in the subspace
$\epsilon^{N-d}L/\epsilon^N L$ of $L/\epsilon^N L$. Moreover, under the
canonical isomorphism $\epsilon^{N-d} L/\epsilon^N L \cong
L/\epsilon^d L $, the subspace $\ker(A_N)$ goes over to the subspace
$\ker(A_d)$. 
\end{lemma}
\begin{proof} 
First we note that
\[
\dim \ker(A_N)=\dim \cok(A_N) =\dim L/AL = d.
\] 
Here we used \eqref{eq.edLc} and $N \ge d$ to see that $AL \supset
\epsilon^N L$. 

 Since
$\det(A)\ne 0$, we may consider the inverse
$A^{-1}
\in M_nF$ of
$A$. From \eqref{eq.edLc} we obtain 
\[
A^{-1}\epsilon^N L \subset \epsilon^{N-d}L \subset L.
\] 
Therefore 
\[
\ker(A_N)=A^{-1}\epsilon^N L/\epsilon^N L \cong 
A^{-1}\epsilon^d L/\epsilon^d L =\ker(A_d)
\]
and 
\[
\ker(A_N)=A^{-1}\epsilon^N L/\epsilon^N L \subset \epsilon^{N-d}
L/\epsilon^N L,
\] 
as desired. 
\end{proof} 

\subsection{Solving the equation $f(x')=y$ by successive
approximations} 
As mentioned before, the next lemma is a variant of  results of Greenberg.
Since the precise statement we need is not stated explicitly in
\cite{greenberg66}, we thought it best to write out in full our adaptation
of Greenberg's arguments. 

\begin{lemma} \label{lem.suc.app}
 Let $x \in L$ and
assume that $g(x) \ne 0$. 
Let
$M$ be an integer such that
$M > d(x)$.  Then 
\[
f(x+\epsilon^{M}L) \supset f(x)+\epsilon^{M+d(x)} L. 
\]
More precisely 
\[
f(x+\epsilon^{M}L) = f(x)+D_x\epsilon^{M} L. 
\]
\end{lemma}

\begin{proof} 
In this proof we abbreviate $d(x)$ to $d$. 
We begin with two observations. The first is
that 
\begin{equation} \label{eq.star}
\epsilon^dL \subset D_xL,
\end{equation}
an instance of \eqref{eq.edLc}. The second is that 
\begin{equation} \label{eq.dagger} 
f(x+h)\equiv f(x)+D_x \cdot h
\mod{\epsilon^{2M}L}
\end{equation}
 for all $h \in \epsilon^{M}L$ (obvious).

It follows from \eqref{eq.star} that $\epsilon^{M+d}L \subset
D_x\epsilon^ML$; therefore the first assertion of the lemma
follows from the second. As for the second, the inclusion 
\[
f(x+\epsilon^M L) \subset f(x)+D_x(\epsilon^ML)
\]
follows from \eqref{eq.dagger} and the fact that 
$\epsilon^{2M}L \subset D_x(\epsilon^M L)$, a consequence of
\eqref{eq.star} and our hypothesis that $M>d$.

It remains only to prove the reverse inclusion, so 
 let $y \in f(x) + D_x(\epsilon^M L)$. We need to find an
element $x' \in x+\epsilon^M L$  such that
$f(x')=y$.
We will obtain $x'$ as the limit of a sequence
$x=x_0,x_1,x_2,\dots$ in $L$ constructed inductively so as to
satisfy the two conditions
\begin{equation}\label{eq.at}
x_i-x_{i-1} \in \epsilon^{M+i-1}L,
\end{equation}
\begin{equation}\label{eq.sharp}
y-f(x_i)\in \epsilon^{M+d+i}L
\end{equation}
for all $i \ge 1$. Some care is needed, because the first step is slightly
different from all the remaining ones. 

We begin by constructing $x_1$. Write $y$ as $f(x)+D_xh$ with $h
\in \epsilon^ML$, and then put $x_1:=x+h$. Clearly \eqref{eq.at}
holds, and \eqref{eq.sharp} follows from \eqref{eq.dagger} and
the hypothesis that $M >d$. 

Now suppose that $i>1$ and that we have already constructed
$x_1,x_2,\dots,x_{i-1}$ satisfying \eqref{eq.at} and
\eqref{eq.sharp}. From \eqref{eq.at} it follows that 
\begin{equation}\label{eq.flat}
x_{i-1}-x \in \epsilon^M L.  
\end{equation}

In particular $g(x_{i-1})\equiv g(x) \mod \epsilon^M \mathcal
O$, and since $M > d$, we conclude that $d(x_{i-1})=d$. Now
applying \eqref{eq.star} to $x_{i-1}$ rather than $x$, we see
that 
\[
\epsilon^{M+d+i-1}L \subset D_{x_{i-1}}\epsilon^{M+i-1}L.
\]
Using this together with \eqref{eq.sharp} for $i-1$, we see that
there exists $h_i \in \epsilon^{M+i-1}L$ such that 
\begin{equation}\label{eq.hiyf}
D_{x_{i-1}}h_i=y-f(x_{i-1}). 
\end{equation}
Put $x_i:=x_{i-1}+h_i$. It is
clear that \eqref{eq.at} holds. It follows from
\eqref{eq.dagger} (with $M$ replaced by $M+i-1$ and $x$ replaced by
$x_{i-1}$) and \eqref{eq.hiyf} that 
\[
f(x_i)\equiv y \mod \epsilon^{2(M+i-1)}L.
\]
This yields \eqref{eq.sharp} since $2(M+i-1) \ge M+d+i$ (use that
$M > d$ and $i >1$). 

It is clear from \eqref{eq.at} that the sequence $x_i$ has a
limit $x'$. It follows from \eqref{eq.flat} that $x' \in
x+\epsilon^M L$. Finally, we see from \eqref{eq.sharp} that
$f(x')=y$. 
\end{proof}

\subsection{Images under $f$ of admissible subsets of $L$}
Admissibility was discussed earlier in the
context of a smooth scheme $X$ over $\mathcal O$. We are now
interested in the case $X=\mathbb A^n$. Thus $X(\mathcal O)=L$,
and we have the notion of admissible subset in $L$. Moreover we
continue with $f:L \to L$ as in the previous subsection. 
\begin{proposition} \label{prop.im.adm}
Let $M$ and $e$ be non-negative integers,  
and let $Z$ be an $M$-admissible subset of $L$ such that
$d(z)\le e$ for all $z \in Z$.  Then the subset 
$f(Z)$ is
$N$-admissible, where $N$ is any positive integer large
enough that $N >2e$ and $N \ge M+e$. 
\end{proposition}
\begin{proof}
This follows immediately from Lemma \ref{lem.suc.app}, applied not to the
integer $M$, but to the integer $N-e$.  
\end{proof}

\subsection{Fibers of $f_N:L/\epsilon^N L \to L/\epsilon^N L$} 
Let $N$ be a positive integer. Our given morphism $\mathbf
f:\mathbb A^n \to \mathbb A^n$ induces a map 
\[
f_N:L/\epsilon^N L \to L/\epsilon^N L
\]
on $\mathcal O/\epsilon^N\mathcal O$-valued points. There is a
commutative square 
\[
\begin{CD}
L @>f>> L \\
@V\pi_N VV @VV\pi_N V \\
L/\epsilon^N L @>f_N>> L/\epsilon^N L
\end{CD}
\]
in which $\pi_N$ is the canonical surjection $L \to L/\epsilon^N
L$. 

We are interested in the fibers of $f_N$. In the linear case of subsection
\ref{sub.lin.case}, as long as $N \ge d$, all fibers of $A_N$ are translates
of
$\ker(A_N)$, a vector space of dimension $d$ that is essentially
independent of $N$ (see Lemma \ref{lem.lin.ker}). Something similar happens
in the non-linear case, with the role of $d$ being played by $d(x)$, but
since $d(x)$ is no longer constant, the situation is necessarily more
complicated. In order to analyze the fiber $f_N^{-1}(\bar y)$ over a point
$\bar y \in L/\epsilon^N L$ we will need to make an assumption (in part (2)
of the next lemma) ensuring that
$d(x) < N/2$ for all points $x$ in  $\pi_N^{-1}(f_N^{-1}(\bar y))$. 

We will make use of the following definitions. For any non-negative integer
$e$ we denote by $L_{\le e}$ the subset of $L$ consisting of all points $x$
for which $d(x) \le e$. The subset $L_{\le e}$ of $L$ is obviously
$(e+1)$-admissible. For $x \in L_{\le e}$ and any integer $N$ such that $N
\ge e$ we may apply Lemma \ref{lem.lin.ker} to the differential $D_x$ to
conclude that $\ker(D_{x,N})$ has dimension $d(x)$, is contained in
$\epsilon^{N-e}L/\epsilon^N L \cong L/\epsilon^e L$, and is independent of
$N$ when viewed as a subspace of $L/\epsilon^e L$. For $x \in L_{\le
e}$ we define 
\[
V_x:= \ker D_{x,e}
\]
and we then have canonical isomorphisms 
\[
V_x\cong \ker D_{x,N}
\] 
for $N \ge e$. 

Now suppose that $M > e$. Since $L_{\le e}$ is $M$-admissible, it is the
preimage of its image $(L/\epsilon^M L)_{\le e}$ in $L/\epsilon^M L$. Since
$D_{x,e}$ depends only on $\pi_e(x) \in L/\epsilon^e L$, so too does the
linear subspace $V_x$  of $L/\epsilon^e L$, so that for
any $z \in (L/\epsilon^M L)_{\le e}$ we get a well-defined linear subspace
$V_z$ of $L/\epsilon^e L$ by putting $V_z:=V_x$ for any $x \in L$ such that
$\pi_M(x)=z$. 

\begin{lemma} \label{lem.fib.s}
Let $e,N$ be non-negative integers satisfying
$N > 2e$. Put $M:=N-e$, noting that $N \ge M >e$. Let $y \in L$ and put
$\bar y:=\pi_N(y) \in L/\epsilon^N L$.
\begin{enumerate}
\item For any $x \in f^{-1}(y)$ the $d(x)$-dimensional affine linear
subspace
\[
A_x:=\pi_N(x)+\epsilon^MV_{\pi_M(x)}
\] 
of $L/\epsilon^N L$ is contained in $f_N^{-1}(\bar y)$. Here
$\epsilon^MV_{\pi_M(x)}$ is the $d(x)$-dimensional linear subspace 
of $\epsilon^ML/\epsilon^N L$ corresponding to $V_{\pi_M(x)}$ under the
canonical isomorphism $\epsilon^M L/\epsilon^N L\cong L/\epsilon^e L$. 
Since the image of $\epsilon^MV_{\pi_M(x)}$ in $L/\epsilon^ML$ is $0$, the
image of $A_x$ in $L/\epsilon^M L$ is the single point $\pi_M(x)$. 
\item If $f_N^{-1}(\bar y)$ is contained in $(L/\epsilon^N L)_{\le e}$, then
$f^{-1}(y)$ is finite and 
\[
f_N^{-1}(\bar y)=\bigcup_{x \in f^{-1}(y)} A_x.
\]
\item If the  composed map
\begin{equation}\label{eq.comp}
f^{-1}(y) \hookrightarrow L \xrightarrow{\pi_{M}} 
L/\epsilon^{M} L
\end{equation}
is injective, then $A_x$, $A_{x'}$ are disjoint whenever $x$, $x'$ are
distinct points in $f^{-1}(y)$.  

\end{enumerate}
\end{lemma}
\begin{proof}
(1) Let 
 $x \in f^{-1}(y)$.  Since $N > 2e$
implies $2M \ge N$, equation \eqref{eq.dagger} tells us
that
\begin{equation}\label{eq.pt1}
f(x+h)\equiv f(x)+D_x \cdot h \mod \epsilon^N L
\end{equation}
for all $h \in \epsilon^{M}L$.  
Using the canonical isomorphism $\epsilon^M L/\epsilon^N L \cong
L/\epsilon^e L$, the restriction of $D_{x,N}$ to $\epsilon^M L/\epsilon^N
L$ becomes identified with $D_{x,e}$, whose kernel is by definition
$V_x=V_{\pi_M(x)}$. 
Therefore $A_x$ 
 is  contained in $f_N^{-1}(\bar y)$.

(2) The morphism $\mathbf f$ is \'etale off the closed subscheme
defined by the vanishing of the Jacobian $g$. It follows that for
any $y' \in L$ there are only finitely many points $x$  in the
fiber
$f^{-1}(y')$ for which $g(x) \ne 0$. Thus,  the
 hypothesis made in (2) ensures the fiber $f^{-1}(y)$ is
indeed finite.

Any element $\bar x' \in f_N^{-1}(\bar y)$ is represented by an element
$x'
\in L$ such that
\begin{equation}\label{eq.fx'y}
f(x')\equiv y \mod \epsilon^N L.
\end{equation}
 Our assumption that $f_N^{-1}(\bar y)$ is contained in $(L/\epsilon^N
L)_{\le e}$ tells us that
$d(x') \le e$. Since $N > 2e$,  Lemma \ref{lem.suc.app} 
 says that
there exists $h \in \epsilon^{M}L$ such that $f(x'-h)=y$. 

Put $x:=x'-h$. Then $x \in f^{-1}(y)$ and we claim that $\bar x' \in A_x$.
Indeed, from
\eqref{eq.pt1} we see that 
\begin{equation}\label{eq.ne.th}
f(x')=f(x+h)\equiv f(x)+D_x\cdot h \mod \epsilon^N L.
\end{equation}
Since $f(x)=y$ and $f(x')\equiv y \mod \epsilon^N L$, 
\eqref{eq.ne.th} implies that $D_x\cdot h \in \epsilon^N L$,
showing that $h$ represents an element in $\ker
(D_{x,N})=\epsilon^MV_{\pi_M(x)}$.  Therefore  $\bar x' \in A_x$.

(3) We have already noted that all points in
$A_x$  have the same image as $x$ in $L/\epsilon^{M}L$. The
injectivity of \eqref{eq.comp} then assures the disjointness of
$A_x$, $A_{x'}$ when $x$, $x'$ are distinct points in
$f^{-1}(y)$. 
\end{proof}

\subsection{The vector bundle $V^d$}\label{sub.vb.vd} 
We retain all the notation of the
previous subsection. In particular (for $M > e$) at each point $z \in
(L/\epsilon^M L)_{\le e}$ we have the vector space $V_z$, whose dimension
depends on $z$. Now fix a non-negative integer $d$ such that $d \le e$ and
consider the $(e+1)$-admissible subset $L_d$ of $L$ consisting of all
elements
$x$ such that $d(x)=d$. Clearly $L_d$ is contained in $L_{\le e}$. Since
$L_d$ is also $M$-admissible, it is the preimage of its image
$(L/\epsilon^M L)_d$ in $L/\epsilon^M L$.  For each point $z \in
(L/\epsilon^M L)_d$ the vector space $V_z$ is $d$-dimensional. 

We claim that we can assemble the vector spaces $V_z$ into a rank $d$
vector bundle $V^d$ over $(L/\epsilon^M L)_d$. Indeed, we just need to
recall the general principle that, given a homomorphism of vector bundles
over a scheme $Y$, and given a locally closed subset $Z$ of $Y$ over which
the homomorphism has constant rank, the pointwise kernels of the
homomorphism assemble into a vector bundle over $Z$. Here we are applying
this general principle to the differential of $f_M:L/\epsilon^ML \to
L/\epsilon^M L$, viewed as a homomorphism from the tangent bundle of
$L/\epsilon^M L$ to itself. 

Note that  the particular choice of $M$ is unimportant, which is why we
have omitted it from the notation. The smallest possible choice is
$e+1$, so we get a vector bundle $V^d$ over $(L/\epsilon^{e+1} L)_d$, and
its pullback by $(L/\epsilon^M L)_d\twoheadrightarrow(L/\epsilon^{e+1}
L)_d$ gives us the vector bundle $V^d$ for $M$.

\section{Proof of Theorem \ref{thm.a.str}}\label{sec.pr.m.th}
Among other
things, we must show that
$\mathbb A(\mathcal O)_s$ is $N$-admissible when $N > 2e$ (abbreviating
$e(w,r)$ to $e$). This follows from Proposition
\ref{prop.im.adm}, once we note that    the
valuation of the Jacobian equals $e$ on $\mathfrak
t_w(\mathcal O)_r$ (see Lemma \ref{lem.val.Jac}), and also that 
$\mathfrak t_w(\mathcal O)_r$ is
$(N-e)$-admissible (use Proposition \ref{prop.t.str}). Here we used the
non-emptiness of our stratum to conclude that for any root $\alpha$ we have
$r(-\alpha)=r(\alpha)$ and hence 
\[
r(\alpha) \le \delta_r/2 \le e.
\]

The rest of the proof is organized as follows. First we will prove the
theorem in the case when $w=1$. Then we will deduce the general case from
this special case. 

\subsection{Proof of Theorem \ref{thm.a.str} when $w=1$} In this subsection
we will always be taking $w=1$. Since we are only interested in non-empty
strata, we only consider functions $r$ on $R$ taking values in the
non-negative integers, and satisfying the property that $R_m:=\{\alpha \in
R: r(\alpha)\ge m\}$ be $\mathbb
Q$-closed for all $m \ge 0$. In this subsection we write simply $\mathbb
A(\mathcal O)_r$ for the stratum in $\mathbb A(\mathcal O)$ obtained as the
image of $\mathfrak t(\mathcal O)_r$. Of course $\mathbb A(\mathcal
O)_{xr}=\mathbb A(\mathcal O)_r$ for all $x \in W$.

The theorem will be a simple consequence of the next two lemmas,
the first involving scaling by $\epsilon^m$, the second involving reduction
to a Levi subgroup.  

In the first of the two lemmas we will need the following
additional notation. For a non-negative integer
$m$ we write
$r+m$ for the function on
$R$ whose value on a root
$\alpha$ is
$r(\alpha)+m$. Also we denote by $d_i$ the
degree of the $i$-th basic invariant $f_i$.

\begin{lemma}\label{lem.red.sca}
Suppose that $\mathbb A(\mathcal O)_r$ is   
locally closed, irreducible,
non-singular of codimension $a$ in $\mathbb A(\mathcal O)$. Then 
$\mathbb A(\mathcal O)_{r+m}$ is   
locally closed, irreducible,
non-singular of codimension $a+m(d_1+\dots +d_n)$ in $\mathbb A(\mathcal
O)$.
\end{lemma}
\begin{proof}
Clearly $\mathfrak t(\mathcal O)_{r+m}=\epsilon^m \mathfrak t(\mathcal
O)_r$. Therefore  $\mathbb A(\mathcal O)_{r+m}=h\mathbb A(\mathcal O)_r$,
where $h$ is the $\mathcal O$-linear map from $L$ to $L$ defined by
$(x_1,\dots,x_n) \mapsto (\epsilon^{md_1}x_1,\dots,\epsilon^{md_n}x_n)$.
Now use Lemma \ref{lem.linear}. 
\end{proof}

The second lemma involves the Levi subgroup $M$ of $G$ containing $T$ whose
root system $R_M$ is equal to the $\mathbb Q$-closed subset 
$R_1=\{\alpha \in R: r(\alpha) \ge 1 \}$ of the root system $R$. We need to
consider $\mathbb A(\mathcal O)$ for both $G$ and $M$, so to avoid
confusion we now write $\mathbb A_G(\mathcal O)$ and $\mathbb A_M(\mathcal
O)$. We write $r_M$ for the function on $R_M$ obtained by restriction from
$r$. 

\begin{lemma}\label{lem.red.lev}
Suppose that $\mathbb A_M(\mathcal O)_{r_M}$ is   
locally closed, irreducible,
non-singular of codimension $a$ in $\mathbb A_M(\mathcal O)$. Then 
$\mathbb A_G(\mathcal O)_{r}$ is  
locally closed, irreducible,
non-singular of codimension $a$ in $\mathbb A_G(\mathcal O)$.
\end{lemma} 
\begin{proof}
We write $W_M$ for the Weyl group of $M$. We will also need the subgroup
$W_M'$ of $W$ defined by 
\[
W_M':=\{w \in W: w(R_M)=R_M \}.
\] 
Note that $W_M$ is a normal subgroup of $W_M'$. Since $\mathbb A_M$ is
obtained by dividing $\mathfrak t$ by the action of the subgroup $W_M$ of
$W$, there is an obvious surjective morphism 
\[
g:\mathbb A_M \to \mathbb A_G,
\]
and this yields a map 
\[
g:\mathbb A_M(\mathcal O) \to \mathbb A_G(\mathcal O).
\]
Note that $W_M'$ acts on $\mathbb A_M$, and that the induced action on
$\mathbb A_M(\mathcal O)$ preserves the fibers of $g$. 

Now define a  polynomial function $Q$ on $\mathfrak t$ by \[
Q=\prod_{\alpha \in R\setminus R_M} \alpha.
\]
Since $Q$ is $W_M$-invariant (even $W_M'$-invariant), we may also regard it
as an element of the ring of regular functions on the affine variety
$\mathbb A_M$. We denote by
$\mathbb A_M^\flat$ the open $k$-subscheme of $\mathbb A_M$ obtained by
removing the locus where $Q$ vanishes. Clearly $\mathbb A_M^\flat$ is stable
under $W_M'$. It follows from the discussion in subsection
\ref{sub.reg.jac} that $Q$ is the square of the Jacobian of $g$ (up to some
non-zero scalar in our base field $k$), and hence that the restriction
$g^\flat$ of
$g$ to
$\mathbb A_M^\flat$ is
\'etale.

Note that $\mathbb A_M^\flat(\mathcal O)$ is
the open, admissible (in fact $1$-admissible) subset of $\mathbb
A_M(\mathcal O)$ consisting of all points $u \in \mathbb A_M(\mathcal
O)$ such that $Q(u)$ is a unit in $\mathcal O$.    
It follows easily from the definitions that  
\begin{equation}\label{eq.disj.N}
(g^\flat)^{-1}(\mathbb A_G(\mathcal O)_r)=\coprod_{x
\in W_M\backslash W_M'}
\bigl(\mathbb A_M(\mathcal O)_{xr_M}\cap \mathbb A_M^\flat(\mathcal
O)\bigr)
\end{equation} 
and 
\begin{equation}\label{eq.gfl.sur}
g^\flat\bigl(\mathbb A_M(\mathcal O)_{r_M} \cap \mathbb A_M^\flat(\mathcal
O)\bigr)=\mathbb A_G(\mathcal O)_r.
\end{equation} 

We claim that for each $x \in W'_M$ the subset 
$\mathbb A_M(\mathcal O)_{xr_M}\cap \mathbb A_M^\flat(\mathcal
O)$  is locally closed, irreducible,
non-singular of codimension $a$ in $\mathbb A_M^\flat(\mathcal O)$. 
Indeed, using the action of $W'_M$ on $\mathbb A_M(\mathcal O)$, we may
assume $x=1$, and then $\mathbb A_M(\mathcal O)_{r_M}\cap \mathbb
A_M^\flat(\mathcal O)$, being open in $\mathbb A_M(\mathcal O)_{r_M}$,
inherits all the stated properties from $\mathbb A_M(\mathcal O)_{r_M}$.
(For irreducibility we need to remark that $\mathbb A_M(\mathcal O)_{r_M}
\cap \mathbb A_M^\flat(\mathcal O)$ is non-empty, because by
\eqref{eq.gfl.sur} it maps onto the non-empty set $\mathbb A_G(\mathcal
O)_r$.)

Next we claim that  $(g^\flat)^{-1}(\mathbb A_G(\mathcal O)_r)$    is
locally closed and  non-singular of codimension $a$ in $\mathbb
A_M^\flat(\mathcal O)$. Using \eqref{eq.disj.N} and the fact that each set
$\mathbb A_M(\mathcal O)_{xr_M} \cap \mathbb A_M^\flat(\mathcal O)$ is
locally closed and non-singular of codimension
$a$, we see that it is enough to show that when $x_1$, $x_2$ are distinct
in $W_M\backslash W_M'$, the closure of $\mathbb A_M(\mathcal O)_{x_1r_M}$
in
$\mathbb A_M(\mathcal O)$ does not meet $\mathbb A_M(\mathcal O)_{x_2r_M}$,
and this follows from Lemma \ref{lem.cl.disj}, applied to $M$ rather than
$G$. 

Now choose a positive integer $N$ so large that $\mathbb A_G(\mathcal O)_r$
is $N$-admissible, and consider the commutative square 
\[
\begin{CD}
\mathbb A_M^\flat(\mathcal O) @>{g^\flat}>> \mathbb A_G(\mathcal O) \\
@VVV @VVV \\
\mathbb A_M^\flat(\mathcal O/\epsilon^N\mathcal O) @>h>> \mathbb
A_G(\mathcal O/\epsilon^N\mathcal O),
\end{CD} 
\] 
in which $h$ is obtained from $g^\flat$ by applying Greenberg's functor and
is therefore \'etale (see 16.1). Since $\mathbb A_G(\mathcal O)_r$ is
$N$-admissible, it is the preimage of its image, call it $\mathbb
A_G(\mathcal O/\epsilon^N\mathcal O)_r$, in $\mathbb A_G(\mathcal
O/\epsilon^N\mathcal O)$. Since $h^{-1}(\mathbb
A_G(\mathcal O/\epsilon^N\mathcal O)_r)$ 
has preimage $(g^\flat)^{-1}(\mathbb A_G(\mathcal O)_r)$ in $\mathbb
A_M^\flat(\mathcal O)$, we conclude that $h^{-1}(\mathbb
A_G(\mathcal O/\epsilon^N\mathcal O)_r)$ is locally closed and non-singular
of codimension $a$ in $\mathbb A_M^\flat(\mathcal O/\epsilon^N\mathcal O)$.
 Since $h$ is \'etale, it is an open
map, and moreover we know from \eqref{eq.gfl.sur} that $\mathbb
A_G(\mathcal O/\epsilon^N\mathcal O)_r$ lies in the image of $h$. It then
follows from 
 Lemmas \ref{lem.top.adm}(1) and \ref{lem.adm.sm} that $\mathbb
A_G(\mathcal O/\epsilon^N\mathcal O)_r$ is locally closed and non-singular. 
Since, again by \eqref{eq.gfl.sur}, $\mathbb
A_G(\mathcal O/\epsilon^N\mathcal O)_r$ is the image of an irreducible
subset of $\mathbb A_M^\flat(\mathcal O/\epsilon^N\mathcal O)$, we conclude
that $\mathbb
A_G(\mathcal O/\epsilon^N\mathcal O)_r$ is irreducible. Since $h$ is
\'etale, the $k$-schemes $\mathbb
A_G(\mathcal O/\epsilon^N\mathcal O)_r$ and $h^{-1}(\mathbb
A_G(\mathcal O/\epsilon^N\mathcal O)_r)$ have the same dimension. 
Moreover $\mathbb
A_G(\mathcal O/\epsilon^N\mathcal O)$ and  $\mathbb A_M^\flat(\mathcal
O/\epsilon^N\mathcal O)$  have the same dimension, so  the
codimension of $\mathbb
A_G(\mathcal O/\epsilon^N\mathcal O)_r$ in
$\mathbb
A_G(\mathcal O/\epsilon^N\mathcal O)$ is the same as that of $h^{-1}(\mathbb
A_G(\mathcal O/\epsilon^N\mathcal O)_r)$ in
$\mathbb A_M^\flat(\mathcal O/\epsilon^N\mathcal O)$, namely $a$. The lemma
is proved. 
\end{proof}
Here is the lemma we needed in the previous proof. It involves two
functions $r,r'$ on $R$ taking values in the non-negative integers. As
usual we assume that all the sets 
\begin{align*}
R_m&=\{\alpha \in R:r(\alpha) \ge m \} \\
R_m'&=\{\alpha \in R: r'(\alpha) \ge m \}
\end{align*}
are $\mathbb Q$-closed. 
\begin{lemma}\label{lem.cl.disj}
Assume that $|R_m|=|R'_m|$ for all $m \ge 0$. If $\mathbb A_G(\mathcal
O)_r$ meets the closure of $\mathbb A_G(\mathcal O)_{r'}$ in $\mathbb
A_G(\mathcal O)$, then
$r'
\in Wr$. 
\end{lemma}
\begin{proof}
Recall from before the integer 
\[
\delta_r=\sum_{\alpha \in R} r(\alpha). 
\]
Our assumption that $|R_m|=|R'_m|$ for all $m \ge 0$ means that $r$,$r'$
have the same set (with multiplicities) of $|R|$ values. In particular
$\delta_r=\delta_{r'}$. 

In this proof we will be using the scalar product 
\[
(r_1,r_2):=\sum_{\alpha \in R}r_1(\alpha)r_2(\alpha),
\]
which is none other than the usual Euclidean inner product on $\mathbb
R^R$. Since, when viewed as vectors in $\mathbb
R^R$, $r$ and $r'$ are permutations of each other, we
have
\[
(r,r)=(r',r').
\]

Pick $N \ge 0$ such that $r(\alpha) \le N$ for all $\alpha \in R$. Define
polynomials $Q_r$, $P_r$ on $\mathfrak t$ by 
\begin{align*}
Q_r&:=\prod_{\alpha \in R}\alpha^{N-r(\alpha)} \\
P_r&:=\sum_{x \in W/W_r}Q_{xr}.
\end{align*}
Since $P_r$ has been defined so as to be $W$-invariant, it can also be
thought of as a regular function on $\mathbb A_G$.

Suppose that $u \in \mathfrak t(\mathcal O)_r$. Then 
\[
\val Q_{xr}(u)=N\delta_r-(xr,r).
\]
Since $(xr,xr)=(r,r)$, the Cauchy-Schwarz inequality implies that $(xr,r) <
(r,r)$ when $xr\ne r$. Therefore 
\[
\val P_{r}(u)=N\delta_r-(r,r).
\]

Now suppose that $u' \in \mathfrak t(\mathcal O)_{r'}$. Then 
\[
\val Q_{xr}(u')=N\delta_{r'}-(xr,r').
\]
Recall that $\delta_r=\delta_{r'}$ and $(r,r)=(r',r')$. The Cauchy-Schwarz
inequality implies that 
\[
\val Q_{xr}(u')\ge N\delta_{r}-(r,r),
\]
with equality only if $xr=r'$; therefore 
\[
\val P_{r}(u')\ge N\delta_{r}-(r,r),
\]
with equality only if $r' \in Wr$. Thus, if $r' \notin Wr$, 
the admissible open subset 
\[\{v \in \mathbb
A_G(\mathcal O): \val P_r(v) \le N\delta_r-(r,r) \}\] 
  of
$\mathbb A_G(\mathcal O)$ contains $\mathbb A_G(\mathcal O)_r$ and is 
disjoint from $\mathbb A_G(\mathcal O)_{r'}$. 
\end{proof}

Now we can prove Theorem \ref{thm.a.str} when $w=1$. One of the assertions
of the theorem is that
$\mathbb A_G(\mathcal O)_r$ has codimension 
\[
D_G(r):=d_G(r)+\frac{1}{2}(\delta_r+c_w)
\]
in $\mathbb A_G(\mathcal O)$, where $d_G(r)$ denotes the codimension of
$\mathfrak t(\mathcal O)_r$ in $\mathfrak t(\mathcal O)$. Since $c_w=0$
when $w=1$, we can write $D_G(r)$ more simply as 
\[
D_G(r)=d_G(r)+\sum_{\alpha \in R^+}r(\alpha).
\] 

We begin by proving parts (1) and (2) of the theorem. 
We do this by induction on $|R|$, the case when $|R|=0$ being
trivial. Now we do the induction step. First suppose that $0$ actually
occurs as a  value of $r$, so that $R_1$ is strictly smaller than $R$.
 Thus the theorem holds for the group $M$ in Lemma
\ref{lem.red.lev} by our inductive hypothesis. Therefore Lemma
\ref{lem.red.lev} implies that $\mathbb A_G(\mathcal O)_{r}$ is   
locally closed, irreducible,
non-singular of codimension $D_M(r_M)$ in $\mathbb A_G(\mathcal O)$. It
remains to check that $D_M(r_M)=D_G(r)$, but this is clear, since $r$
vanishes on roots of $G$ that are not roots of $M$. (Note that $\mathfrak
t(\mathcal O)_r$ is open in $\mathfrak t(\mathcal O)_{r_M}$, so that
$d_G(r)=d_M(r_M)$.)  

Now consider the general case. 
Let $m$ be the smallest integer which actually occurs as a value
of $r$. Then $r$ can be written as $r'+m$, and $0$ occurs as a value of
$r'$, so that the theorem holds for $r'$ by what we have already proved.
Therefore Lemma
\ref{lem.red.sca} implies that $\mathbb A_G(\mathcal O)_{r}$ is   
locally closed, irreducible,
non-singular of codimension $D_G(r')+m(d_1+\dots+d_n)$ in $\mathbb
A_G(\mathcal O)$. It
follows easily from \eqref{eq.rt.d} that $D_G(r')+m(d_1+\dots+d_n)=D_G(r)$,
and the proof of parts (1) and (2) of the theorem is now complete. 

It remains to prove part (3) of the theorem, which asserts that $\mathfrak
t(\mathcal O)_r$ is smooth over $\mathbb A(\mathcal O)_r$.  Since
$\mathfrak t(\mathcal O)_r$ and $\mathbb A(\mathcal O)_r$ are
non-singular, it suffices (see Lemma \ref{lem.sm.dif}) to check that for
each
$u
\in \mathfrak t(\mathcal O)_r$ the
differential $df_u$ of
$f$ maps 
$T_{\mathfrak
t(\mathcal O)_r,u}$ onto $T_{\mathbb A(\mathcal O)_r,f(u)}$. For this we
just need to show that $T_{\mathbb A(\mathcal O)_r,f(u)}$ and the 
image under $df_u$ of
$T_{\mathfrak t(\mathcal O)_r,u}$ have the same codimension in $T_{\mathbb
A(\mathcal O),f(u)}$. From part (2) of the theorem we know that the
codimension of $T_{\mathbb A(\mathcal O)_r,f(u)}$ in
$T_{\mathbb A(\mathcal O),f(u)}$ is 
\begin{equation}\label{eq.codm.fo}
d_G(r)+\sum_{\alpha \in R^+}r(\alpha).
\end{equation}
Since the valuation of the determinant of $df_u$ is $\sum_{\alpha \in
R^+}r(\alpha)$ (see \eqref{eq.jac.pro}), and the codimension of the tangent
space $T_{\mathfrak
t(\mathcal O)_r,u}$ in $T_{\mathfrak t(\mathcal O),u}$ is equal to 
$d_G(r)$, we conclude that the codimension of the image under $df_u$ of
$T_{\mathfrak t(\mathcal O)_r,u}$ is also equal to the expression
\eqref{eq.codm.fo}, and we are done.

\subsection{Proof of Theorem \ref{thm.a.str} in general}
Now we prove the theorem in the general case. So consider a pair $(w,r)$
such that $\mathbb A(\mathcal O)_{w,r}$ is
non-empty, and let $s$ denote the $W$-orbit of $(w,r)$. 

 Here is the idea of the proof.  
 As usual we denote by 
$E$  the field $F_{l}$, with $l=o(w)$. We once again  denote by  
$r_E$  the   integer valued function on $R$ obtained by multiplying
$r$ by
$l$. From the special case of the theorem that we have already proved
(applied to $E$ rather than $F$), we understand
$\mathfrak t(\mathcal O_E)_{r_E}$ and $\mathbb A(\mathcal O_E)_{r_E}$, and
we are going to deduce the theorem in general by taking fixed points of
suitable  automorphisms of order $l$. 

In the case of $\mathfrak t(\mathcal O_E)_{r_E}$, we consider the
automorphism $u \mapsto w\tau_E(u)$. The fixed point set of this action
is of course  $\mathfrak t_w(\mathcal O)_r$.  In the case of $\mathbb
A(\mathcal O_E)_{r_E}$, we consider the action of $\tau_E$. The fixed point
set of this action contains $\mathbb A(\mathcal O)_{w,r}$. These two
automorphisms give us actions of the cyclic group $\mathbb Z/l\mathbb Z$,
and the map 
\begin{equation}\label{eq.pre.fp}
f:\mathfrak t(\mathcal O_E)_{r_E} \to \mathbb
A(\mathcal O_E)_{r_E}
\end{equation} 
is $\mathbb Z/l\mathbb Z$-equivariant. Taking fixed points under $\mathbb
Z/l\mathbb Z$, we get 
\[
\mathfrak t_w(\mathcal O)_r \xrightarrow{f_w} \mathbb A(\mathcal O)_{s}
\hookrightarrow \mathbb
A(\mathcal O_E)_{r_E}^{\mathbb Z/l\mathbb Z}.
\]

Since $l$ is invertible in $k$, taking fixed points under  $\mathbb
Z/l\mathbb Z$ preserves non-singularity and smoothness, as 
is discussed in appendix \ref{sec.fp.gx}. This will end up giving us a good
handle on $\mathbb A(\mathcal O)_s$.  Unfortunately we cannot apply appendix
\ref{sec.fp.gx} directly to
\eqref{eq.pre.fp}, since we need to be dealing with schemes of finite
type over
$k$.  To achieve this we use that all the sets under consideration are
$N$-admissible for sufficiently large $N$. 

More precisely there are four admissible sets under consideration. We begin
by choosing $M$ large enough that $r(\alpha) < M$ for all $\alpha \in R$.
This guarantees (see Proposition \ref{prop.t.str}) that $\mathfrak
t_w(\mathcal O)_r$ and $\mathfrak
t(\mathcal O_E)_{r_E}$ are $M$-admissible. Increasing $M$ as need be, 
we may also assume that $\mathbb A(\mathcal O)_s$ and 
$\mathbb A(\mathcal O_E)_{r_E}$ are $M$-admissible. Now let $N$ be any
integer such that $N \ge M$.

Thus, now letting $\mathfrak
t_w(\mathcal O/\epsilon^N\mathcal O)_r$ denote the image of $\mathfrak
t_w(\mathcal O)_r$ under $\mathfrak t_w(\mathcal O) \twoheadrightarrow
\mathfrak t_w(\mathcal O/\epsilon^N\mathcal O)$, the set $\mathfrak
t_w(\mathcal O)_r$ is the preimage of $\mathfrak t_w(\mathcal
O/\epsilon^N\mathcal O)_r$. Similarly,  letting $\mathfrak t(\mathcal
O_E/\epsilon^N\mathcal O_E)_{r_E}$ denote the image of $\mathfrak
t(\mathcal O_E)_{r_E}$ under $\mathfrak t(\mathcal O_E) \twoheadrightarrow
\mathfrak t(\mathcal O_E/\epsilon^N\mathcal O_E)$, the set $\mathfrak
t(\mathcal O_E)_{r_E}$ is the preimage of $\mathfrak t(\mathcal
O_E/\epsilon^N\mathcal O_E)_{r_E}$. 

In addition $\mathbb A(\mathcal O)_s$ is
the preimage of its image $\mathbb A(\mathcal O/\epsilon^N\mathcal O)_s$ in
$\mathbb A(\mathcal O/\epsilon^N\mathcal O)$, and similarly  
$\mathbb A(\mathcal O_E)_{r_E}$ is
the preimage of its image $\mathbb A(\mathcal O_E/\epsilon^N\mathcal
O_E)_{r_E}$ in
$\mathbb A(\mathcal O_E/\epsilon^N\mathcal O_E)$. With all this notation
in place, we can now finish the proof. 

As noted in subsection \ref{sub.gre.re.sc}, there are two different ways to
use Greenberg's functor to regard $\mathfrak t(\mathcal
O_E/\epsilon^N\mathcal O_E)$ as the set of $k$-points of a $k$-scheme. One
is to apply Greenberg's functor directly to $\mathfrak
t$, but working with $\mathcal O_E$ rather than $\mathcal O$. The other is
to apply restriction of scalars $R_{\mathcal O_E/\mathcal O}$ to $\mathfrak
t$ and then use Greenberg's functor for $\mathcal O$. Fortunately,
\ref{sub.gre.re.sc} assures us that the two methods give the same result,
so we will be free to use whichever interpretation is  most convenient at a
given moment. The same remarks apply to $\mathbb A(\mathcal
O_E/\epsilon^N\mathcal O_E)$.

Consider the commutative square
\[
\begin{CD}
\mathfrak t(\mathcal
O_E/\epsilon^N\mathcal O_E)_{r_E} @>>> \mathbb A(\mathcal O_E
/\epsilon^N \mathcal O_E)_{r_E} \\
@VVV @VVV \\
\mathfrak t(\mathcal
O_E/\epsilon^N\mathcal O_E) @>>> \mathbb A(\mathcal O_E
/\epsilon^N \mathcal O_E).
\end{CD} 
\]
From Proposition \ref{prop.t.str} and the special case of the theorem that
has already been proved, we know that the vertical arrows are locally closed
immersions, that the top horizontal arrow is smooth, and that all four
corners of the square are non-singular.  

Recall from \ref{sub.tw.fp} the $\mathbb Z/l\mathbb Z$-action on
$R_{\mathcal O_E/\mathcal O}\mathfrak t$  whose fixed point
scheme is $\mathfrak t_w$. From it we get an action of $\mathbb Z/l\mathbb
Z$ on the $k$-scheme 
\[
(R_{\mathcal O_E/\mathcal
O}\mathfrak t)(\mathcal O/\epsilon^N\mathcal O)=\mathfrak t(\mathcal
O_E/\epsilon^N\mathcal O_E).
\] 
It follows from  Proposition
\ref{prop.ne.gen}(3) that our action preserves $\mathfrak t(\mathcal
O_E/\epsilon^N\mathcal O_E)_{r_E}$ set-theoretically, hence 
scheme-theoretically as well, since we are using the induced reduced
subscheme structure.
Similarly, $\mathbb Z/l\mathbb Z$ acts on $\mathbb A(\mathcal O_E
/\epsilon^N \mathcal O_E)$, preserving the locally closed subscheme $\mathbb A(\mathcal O_E
/\epsilon^N \mathcal O_E)_{r_E}$. 

Now we take fixed points under $\mathbb Z/l\mathbb Z$ everywhere in the 
commutative square we are considering. Bearing in mind that taking fixed
points commutes with Greenberg's functor (see subsection
\ref{sub.gre.fps}), we obtain the commutative square 
\[
\begin{CD}
(\mathfrak t(\mathcal
O_E/\epsilon^N\mathcal O_E)_{r_E})^{\mathbb Z/l\mathbb Z} @>>> (\mathbb
A(\mathcal O_E /\epsilon^N \mathcal O_E)_{r_E} )^{\mathbb Z/l\mathbb Z}\\
@VVV @VVV \\
\mathfrak t_w(\mathcal
O/\epsilon^N\mathcal O) @>>> \mathbb A(\mathcal O
/\epsilon^N \mathcal O).
\end{CD} 
\]
Since taking $\mathbb Z/l\mathbb Z$-fixed points preserves immersions,
non-singularity, and smoothness (see Lemma \ref{lem.Gsmth}), we conclude
that all four corners of our square are non-singular, that the top
horizontal arrow is smooth, and that the two vertical arrows are locally
closed immersions. 

Using that $N > r(\alpha)$ for all $\alpha \in R$, one sees easily that
$(\mathfrak t(\mathcal
O_E/\epsilon^N\mathcal O_E)_{r_E})^{\mathbb Z/l\mathbb Z}$ coincides with
$\mathfrak t_w(\mathcal O/\epsilon^N\mathcal O)_r$ set-theoretically. Since
both are non-singular schemes, hence reduced, they actually coincide as
subschemes. 

The image of $(\mathfrak t(\mathcal
O_E/\epsilon^N\mathcal O_E)_{r_E})^{\mathbb Z/l\mathbb Z}=\mathfrak
t_w(\mathcal O/\epsilon^N\mathcal O)_r$ in $(\mathbb A(\mathcal
O_E/\epsilon^N\mathcal O_E)_{r_E})^{\mathbb Z/l\mathbb Z}$ is open (since
the top horizontal arrow is smooth) and its further image in $\mathbb
A(\mathcal O/\epsilon^N\mathcal O)$, namely $\mathbb A(\mathcal
O/\epsilon^N\mathcal O)_{s}$, is therefore locally closed in $\mathbb
A(\mathcal O/\epsilon^N\mathcal O)$. At the same time we see that 
$\mathbb A(\mathcal
O/\epsilon^N\mathcal O)_{s}$ is non-singular and that 
\begin{equation}\label{eq.nice.smN}
\mathfrak t_w(\mathcal O/\epsilon^N\mathcal O)_r \twoheadrightarrow \mathbb
A(\mathcal O/\epsilon^N\mathcal O)_{s}
\end{equation}
is smooth. Since we have proved that \eqref{eq.nice.smN} is smooth for all
$N\ge M$, we conclude that  
\begin{equation}\label{eq.nice.sm}
\mathfrak t_w(\mathcal O)_r \twoheadrightarrow \mathbb A(\mathcal O)_{s}
\end{equation} 
is smooth. 
Since
$\mathfrak t_w(\mathcal O/\epsilon^N\mathcal O)_r$ is irreducible, so too
is  $\mathbb A(\mathcal O/\epsilon^N\mathcal O)_{s}$. 

At this point we have proved all parts of the theorem except for the
statement concerning the codimension of $\mathbb A(\mathcal O)_s$. For this
we use tangent spaces (which we are free to use since we now know
that the admissible subsets $\mathfrak t_w(\mathcal O)_r$ and  $\mathbb
A(\mathcal O)_{s}$ are locally closed and
non-singular). Choose some point $u \in \mathfrak t_w(\mathcal O)_r$ and
let $c$ denote its image in $\mathbb A(\mathcal O)_s$. The codimension of
$\mathbb A(\mathcal O)_s$ in
$\mathbb A(\mathcal O)$ is the same as that of the tangent space to
$\mathbb A(\mathcal O)_s$ at $c$ in the tangent space to $\mathbb
A(\mathcal O)$ at $c$. Now, since \eqref{eq.nice.sm} is smooth, Lemma
\ref{lem.sm.dif} tells us that the tangent space to $\mathfrak t_w(\mathcal
O)_r$ at
$u$ maps onto the tangent space   to $\mathbb A(\mathcal O)_s$ at $c$. 

We conclude that the codimension of $\mathbb A(\mathcal O)_s$ in
$\mathbb A(\mathcal O)$ is the sum of the codimension of $\mathfrak
t_w(\mathcal O)_r$ in $\mathfrak t_w(\mathcal O)$ (a number we have denoted
by
$d(w,r)$) and the valuation of the Jacobian of $f_w$ at the point $u$
(which by Lemma \ref{lem.val.Jac} we know to be equal to
$e(w,r)=(\delta_r+c_w)/2$). This finally finishes the proof of the theorem. 

\section{Proof of Theorem \ref{thm.HH}}\label{sec.pf.H} 
The idea of the proof is simple enough. We will check that Lemma
\ref{lem.fib.s} applies to our situation, concluding that each
fiber of  the morphism $\mathfrak t_w(\mathcal O/
\epsilon^N\mathcal O)_r \to \mathbb A(\mathcal
O/\epsilon^N\mathcal O)_s$
is a disjoint union of affine spaces of dimension $e$. These affine spaces
are permuted simply  transitively by $W_{w,r}$, and we have already proved
that the morphism is smooth. This makes it plausible that the theorem is
true, but we must construct the rank $e$ vector bundle $\tilde V$ and check
that the morphism really is a torsor for $H=W_{w,r}\times \tilde V$.  

Let $c \in \mathbb A(\mathcal
O)_s$ and let $\bar c$ denote the image of $c$ in $\mathbb A(\mathcal
O/\epsilon^N\mathcal O)_s$. 
 Lemma \ref{lem.fib.s} will give us information about the
fiber $Z$ of the morphism 
\[
f_{w,N}:\mathfrak t_w(\mathcal O/
\epsilon^N\mathcal O) \to \mathbb A(\mathcal
O/\epsilon^N\mathcal O)
\]
over the point $\bar c$. 
 
We need to check that the hypotheses of the lemma are verified. As in that
lemma we will use $L$ to denote $\mathcal O^n=\mathbb A(\mathcal O)$. We
are assuming that $N>2e$, so we just need to verify the assumptions made in
(2) and (3) of the lemma. 

By Lemma \ref{lem.inv.i.as} we have  
\begin{equation}\label{eq.d.u.f}
 f_w^{-1}(\mathbb A(\mathcal
O)_s)=\coprod_{x \in W_w/W_{w,r}}\mathfrak
t_w(\mathcal O)_{xr}.
\end{equation}
It then follows from  Lemma
\ref{lem.val.Jac} that $\val \det df_w$ takes the
constant value $e$ on $f_w^{-1}(\mathbb A(\mathcal
O)_s)$. Since $\mathbb A(\mathcal O)_s$ is $N$-admissible, the preimage
of $Z$ in $\mathfrak t_w(\mathcal O)$ is contained in $f_w^{-1}(\mathbb
A(\mathcal O)_s)$, and therefore  $\val \det df_w$ takes the
constant value $e$ on that preimage,  showing that the assumption
about 
$f_{w,N}^{-1}(\bar c)$ made in (2) of Lemma \ref{lem.fib.s} does hold. 

As for the assumption on the fiber $f_w^{-1}(c)$ made in (3) of Lemma
\ref{lem.fib.s}, we first recall (see \ref{sub.fib.TWA}) that the group
$W_w$ acts simply transitively on this fiber.
 Next, recall from the first paragraph of section \ref{sec.pr.m.th} that 
\begin{equation}\label{eq.r.ineq}
r(\alpha) \le e < N-e \qquad \forall \, \alpha \in R
\end{equation}
and hence that $\mathfrak t_w(\mathcal O)_r$ is $(N-e)$-admissible. 
Since (see 
\ref{sub.free} and use \eqref{eq.r.ineq})  
$W_w$ acts freely on the image of this fiber in $\mathfrak t_w(\mathcal O/
\epsilon^{N-e}\mathcal O)$, we conclude that the fiber injects into 
$\mathfrak t_w(\mathcal O/
\epsilon^{N-e}\mathcal O)$, as desired. 

The lemma then describes the fiber $Z=f_{w,N}^{-1}(\bar c)$  as a disjoint
union of affine spaces
$A_u$, one for each $u \in f_w^{-1}(c)$. However, we are really interested
in the fiber $g^{-1}(\bar c)$ of the morphism 
\[
g:\mathfrak t_w(\mathcal O/
\epsilon^N\mathcal O)_r \to \mathbb A(\mathcal
O/\epsilon^N\mathcal O)_s
\]
obtained by restriction from $f_{w,N}$. 
For each $u \in f_w^{-1}(c)$ there exists (by \eqref{eq.d.u.f}) $x \in W_w$
such that $u \in \mathfrak t_w(\mathcal O)_{xr}$. We noted in part (1) of 
Lemma
\ref{lem.fib.s} that  all the points in $A_u$ have the same image as $u$
in $\mathfrak t_w(\mathcal O/\epsilon^{N-e}\mathcal O)$. Since $\mathfrak
t_w(\mathcal O)_{xr}$ is $(N-e)$-admissible, it follows that $A_u \subset
\mathfrak t_w(\mathcal O/\epsilon^N\mathcal O)_{xr}$. 
Therefore 
\begin{equation}\label{eq.d.u.5}
g^{-1}(\bar c) =Z\cap \mathfrak t_w(\mathcal O/\epsilon^N\mathcal
O)_r=\coprod_{u} A_u,
\end{equation} 
where the index set for the disjoint union is $f_w^{-1}(c)\cap\mathfrak
t_w(\mathcal O)_r$. From Lemma \ref{lem.inv.i.as} we know that $W_{w,r}$
acts simply transitively on $f_w^{-1}(c)\cap\mathfrak
t_w(\mathcal O)_r$. Thus the natural action of
$W_{w,r}$ on
$g^{-1}(\bar c)$ permutes simply transitively the $e$-dimensional affine
spaces
$A_u$ appearing in the disjoint union  \eqref{eq.d.u.5}.  

These affine spaces arise as orbits of  translation actions of certain
vector spaces described in Lemma \ref{lem.fib.s}. We are going to use the
discussion in
\ref{sub.vb.vd} to assemble these vector spaces into a vector bundle.
Eventually we will arrive at the vector bundle
$\tilde V$, but we must begin with the one (over a different base space)
that is provided by \ref{sub.vb.vd}. 

Put $M:=N-e$. We have already noted  that $\mathfrak
t_w(\mathcal O)_r$ is
$M$-admissible. Therefore the obvious surjection
\[
\pi:\mathfrak t_w(\mathcal O/
\epsilon^N\mathcal O)_r \to \mathfrak t_w(\mathcal O/
\epsilon^M\mathcal O)_r
\]
 is an affine space bundle, more precisely, a torsor (actually trivial,
not that it matters) under the vector group 
\begin{equation}\label{eq.ve.gr}
\ker[\mathfrak t_w(\mathcal O/
\epsilon^N\mathcal O) \to \mathfrak t_w(\mathcal O/
\epsilon^M\mathcal O)]=\mathfrak t_w(\mathcal O/
\epsilon^e\mathcal O)
\end{equation}
(the identification being made using multiplication by $\epsilon^M$). 

Since $\val\det df_w$ takes the constant value $e$ on $\mathfrak
t_w(\mathcal O)_r$, subsection \ref{sub.vb.vd} provides us with a rank $e$
vector bundle $V$  over $\mathfrak
t_w(\mathcal O/\epsilon^M\mathcal O)_r$, obtained by restriction from the
vector bundle $V^e$ of \ref{sub.vb.vd}. In fact
$V$ is a subbundle of the constant vector bundle over $\mathfrak
t_w(\mathcal O/\epsilon^M\mathcal O)_r$ with fiber \eqref{eq.ve.gr}. It is
clear from its definition  that $V$ is $W_{w,r}$-equivariant with respect
to the natural action of
$W_{w,r}$ on $\mathfrak t_w(\mathcal O/\epsilon^M\mathcal O)_r$.

The vector bundle $V$ acts by translations on the affine space bundle 
$\mathfrak t_w(\mathcal O/
\epsilon^N\mathcal O)_r$ over $\mathfrak t_w(\mathcal O/
\epsilon^M\mathcal O)_r$, and we may divide out by its action, obtaining a
factorization 
\[
\mathfrak t_w(\mathcal O/
\epsilon^N\mathcal O)_r \xrightarrow{\rho} \mathfrak t_w(\mathcal O/
\epsilon^N\mathcal O)_r/V \xrightarrow{\eta}\mathfrak t_w(\mathcal O/
\epsilon^M\mathcal O)_r
\]
of $\pi$, in which $\rho$, $\eta$ are both affine space bundles. More
precisely $\rho$ is a torsor for $\eta^*V$, and $\eta$ is a torsor for the
vector bundle obtained by taking the quotient of the constant vector bundle
$\mathfrak t_w(\mathcal O/
\epsilon^e\mathcal O)$ by its subbundle $V$. 

By Lemma \ref{lem.fib.s} the morphism $g$ is constant on the fibers of the
bundle
$\rho$. By faithfully flat descent we see that $g$ factors uniquely as 
\[
\mathfrak t_w(\mathcal O/
\epsilon^N\mathcal O)_r \xrightarrow{\rho} \mathfrak t_w(\mathcal O/
\epsilon^N\mathcal O)_r/V \xrightarrow{h}\mathbb A(\mathcal O/
\epsilon^N\mathcal O)_s.
\]
(To apply descent theory we just need to check the equality of two
morphisms $B \to \mathbb A(\mathcal O/
\epsilon^N\mathcal O)_s$, where $B$ denotes the fiber product of 
$\mathfrak t_w(\mathcal O/
\epsilon^N\mathcal O)_r$ with itself over $\mathfrak t_w(\mathcal O/
\epsilon^N\mathcal O)_r/V$. Now $B$, being itself an affine space bundle
over the reduced scheme $\mathfrak t_w(\mathcal O/
\epsilon^M\mathcal O)_r$, is also reduced, so that the equality of our two
morphisms $B \to \mathbb A(\mathcal O/
\epsilon^N\mathcal O)_s$ follows from the obvious fact that they coincide
on $k$-points.)

Now $g$ is smooth (by Theorem \ref{thm.a.str}) and so is $\rho$; therefore
$h$ is smooth as well. The $W_{w,r}$-equivariance of $V$ ensures that the
action of $W_{w,r}$ on $\mathfrak t_w(\mathcal O/
\epsilon^N\mathcal O)_r$ descends to an action on $\mathfrak t_w(\mathcal O/
\epsilon^N\mathcal O)_r/V$ over $\mathbb A(\mathcal O/
\epsilon^N\mathcal O)_s$, and Lemma \ref{lem.fib.s} tells us that $W_{w,r}$
acts simply transitively on the fibers of $h$. This means that $h$ is in
fact \'etale, and hence that $\mathfrak t_w(\mathcal O/
\epsilon^N\mathcal O)_r/V$ is a $W_{w,r}$-torsor over $\mathbb A(\mathcal O/
\epsilon^N\mathcal O)_s$. 

The pullback $\eta^*V$ is a $W_{w,r}$-equivariant vector bundle over 
$\mathfrak t_w(\mathcal O/
\epsilon^N\mathcal O)_r/V$. Since $h$ is a $W_{w,r}$-torsor, $\eta^*V$
descends to a vector bundle $\tilde V$ on 
 $\mathbb A(\mathcal O/
\epsilon^N\mathcal O)_s$, and we see from the factorization $g=h\rho$ that
$g$ is a $(W_{w,r}\times \tilde V$)-torsor. 
The proof is now complete.

\section{Appendix. Technical lemmas related to admissibility}
In this appendix we verify some lemmas  needed to back up the
statements we made in section \ref{sec.adm.sub} concerning admissible
subsets of
$X(\mathcal O)$. 

\subsection{Elementary facts about open mappings} 
\begin{lemma}\label{lem.adm.cl}
Let $f:Y \to X$ be a continuous map of topological spaces. Then the
following three conditions are equivalent: 
\begin{enumerate}
\item $f$ is an open mapping.
\item For every closed subset $Z \subset Y$ the set $\{ x \in X: f^{-1}(x)
\subset Z \}$ is closed in $X$. 
\item For every subset $S \subset X$ we have
$\overline{f^{-1}(S)}=f^{-1}(\overline{S})$. Here the overlines indicate
closures. 
\end{enumerate}  
\end{lemma}
\begin{proof}
(1) holds iff $f(U)$ is open for every open  $U \subset Y$. Phrasing this
in complementary terms,  (1) holds iff $f(Z^c)^c$ is
closed for every closed subset $Z \subset Y$, where the superscript $c$
indicates  complement. Since $f(Z^c)^c=\{ x \in X: f^{-1}(x)
\subset Z \}$, we see that (1) is equivalent to (2). 

Now consider (3).  Since $f$ is continuous,  
$f^{-1}\overline S$ is a closed subset containing $f^{-1}(S)$.
Therefore (3) holds iff for every $S \subset X$ and every closed $Z \subset
Y$ we have the implication $Z \supset f^{-1}S \Longrightarrow Z \supset
 f^{-1}\overline{S}$. This last implication can be rewritten as $S \subset 
\{ x \in X: f^{-1}(x)
\subset Z \} \Longrightarrow \overline{S} \subset \{ x \in X: f^{-1}(x)
\subset Z \}$, which makes it clear that the implication holds for all $S$
iff
$\{ x
\in X: f^{-1}(x)
\subset Z \}$ is closed. Therefore (3) is equivalent to (2). 
\end{proof}

\begin{lemma}\label{lem.top.adm}
Let $f:Y\to X$ be a continuous, open, surjective map of topological spaces,
and let $S$ be a subset of $X$. Then 
\begin{enumerate}
\item The set $S$ is open \textup{(}respectively, closed, locally
closed\textup{)} in
$X$  iff
$f^{-1}S$ is open \textup{(}respectively, closed, locally closed\textup{)}
in
$Y$. 
\item Assume further that each fiber of $f$ is an irreducible topological
space. Then $S$ is irreducible iff $f^{-1}S$ is irreducible.  
\end{enumerate}
\end{lemma}
\begin{proof}
(1) Everything here is well-known (and obvious) except possibly the fact
that if $f^{-1}S$ is locally closed, then $S$ is locally closed. So suppose
that $f^{-1}S$ is locally closed, which means  that $f^{-1}S$ is open
in its closure. Using (3) in the previous lemma, we see that $f^{-1}S$
is open in 
$f^{-1}\overline{S}$. Since the map $f^{-1}\overline{S} \to \overline{S}$
(obtained by restriction from $f$) is obviously open, we conclude that
$ff^{-1}S=S$ is open in $\overline S$, which means that
$S$ is locally closed. 

(2) $(\Longleftarrow)$ Clear. $(\Longrightarrow)$ Now
assume all fibers of $f$ are irreducible, and assume further that $S$ is
irreducible. We must show that $f^{-1}S$ is irreducible, so suppose that
$Y_1$, $Y_2$ are closed subsets of $Y$ such that $f^{-1}S \subset Y_1 \cup
Y_2$. Put $X_i:=\{x \in X : f^{-1}(x) \subset Y_i\}$ for $i=1,2$. We know
from the previous lemma that $X_1$ and $X_2$ are closed in $X$, and using
the irreducibility of the fibers of $f$, we see that $S \subset X_1 \cup
X_2$. Since $S$ is irreducible, it follows that $S \subset X_1$ or $S
\subset X_2$. Therefore $f^{-1}S \subset
f^{-1}X_1 \subset Y_1$ or $f^{-1}S \subset
f^{-1}X_2 \subset Y_2$, as desired. 
\end{proof}

\subsection{Lemma on smooth morphisms}
\begin{lemma}\label{lem.adm.sm}
Let $X$, $Y$ be schemes locally of finite type over a  noetherian base
scheme
$S$. Let $f:Y \to X$ be a smooth  $S$-morphism. Let $X'$ be a
locally closed subset of $X$, let $Y'$ denote the locally closed subset
$f^{-1}X'$ of $Y$, and equip both $X'$ and $Y'$ with their induced reduced
subscheme structures. Then the natural morphism $Y' \to Y \times_X X'$ is
an isomorphism, and $Y'$ is smooth over
$X'$. If in addition $X' \subset fY$, then 
$X'$ is smooth over $S$ if and only if $Y'$ is smooth over $S$. 
\end{lemma} 
\begin{proof}
First note that $Y \times_X X'$ is a subscheme of $Y$ with the same
underlying topological space as $Y'$. Moreover $Y \times_X X'$ is smooth
over the reduced scheme $X'$, and therefore (EGA IV (17.5.7)) $Y \times_X
X'$ is reduced, which implies that $Y'=Y\times_X X'$ as closed subschemes.
In particular the morphism  $Y' \to X'$ is smooth. If in addition $X'
\subset fY$, then $Y' \to X'$ is also surjective, and it then follows from
EGA IV (17.11.1) that
$Y'$ is smooth over
$S$ if and only if $X'$ is smooth over $S$. 
\end{proof}

\section{Appendix: Some results of Steinberg} 
In \cite{steinberg75} Steinberg proves a number of delicate results on the
behavior of conjugacy classes in the Lie algebra of 
$G$ when the  characteristic of the base field $k$ is not a torsion prime
for
$G$. In this paper we are operating under the very strong hypothesis that
$|W|$ be invertible in $k$, and this makes life rather simple. Nevertheless
it is convenient to obtain what we need as an easy consequence of
\cite{steinberg75}. 

\subsection{Set-up} Let $S$ be any subset of $\mathfrak t$. Define a subset
$R_S$ of our root system $R$ by 
\[
R_S:=\{\alpha \in R: \alpha(u)=0 \quad \forall \, u \in S \}.
\]
Even with no assumption on the characteristic of $k$, it is clear that
$R_S$ is $\mathbb Z$-closed, in the sense that if $\alpha \in R$ lies in
the $\mathbb Z$-linear span of $R_S$ in $X^*(T)$, then $\alpha \in R_S$. In
particular
$R_S$ is a root system in its own right, whose Weyl group we denote by
$W(R_S)$, a subgroup of $W$ which clearly lies inside the subgroup 
\[
W_S:=\{ w \in W: w(u)=u \quad \forall \, u \in S \}. 
\]

\begin{proposition}\label{prop.stein}
Assume, as usual, that $|W|$ be invertible in $k$. Then 
\begin{enumerate}
\item The subgroups $W_S$ and $W(R_S)$ coincide.
\item The subset $R_S$ is $\mathbb Q$-closed, in the sense that 
if $\alpha \in R$ lies in
the $\mathbb Q$-linear span of $R_S$ in $X^*(T)$, then $\alpha \in R_S$. 
\item There is a Levi subgroup $M \supset T$ in $G$ whose root system $R_M$
coincides with $R_S$. 
\end{enumerate}
\end{proposition}
\begin{proof}
(1) This follows immediately from Corollary 2.8, Lemma 3.7, Corollary 3.11
and Theorem 3.14 in Steinberg's article \cite{steinberg75}. 

(2) Let $L(R)$ (respectively, $L(R_S)$) denote the $\mathbb Z$-linear span
of $R$ (respectively, $R_S$) in $X^*(T)$.  
Similarly, let $L(R^\vee)$ (respectively, $L(R^\vee_S)$) denote the $\mathbb
Z$-linear span of $R^\vee$ (respectively, $R^\vee_S$) in $X_*(T)$. 
Using  a suitably normalized $W$-invariant $\mathbb Z$-valued symmetric
bilinear form on $L(R^\vee)$, we obtain a $W$-equivariant embedding 
\[
\varphi: L(R^\vee) \to L(R) 
\]
such that for every $\alpha \in R$ there exists a positive integer
$d_\alpha$ dividing $|W|$ (hence invertible in $k$) such that
$\varphi(\alpha^\vee)=d_\alpha \alpha$. (In fact we can arrange that
$d_\alpha$ is always $1$, $2$, or $3$, with $3$ occurring only when one of
the irreducible components of $R$ is of type $G_2$.)

Now suppose that $\alpha \in R$
lies in the $\mathbb Q$-linear span of $R_S$ in $X^*(T)$. 
We must show that
$\alpha \in R_S$. Since $\varphi$ becomes an isomorphism after tensoring
with $\mathbb Q$, it is also true that $\alpha^\vee$ lies in the $\mathbb
Q$-linear span of $R_S^\vee$ in $X_*(T)$. Therefore  
the class of $\alpha^\vee$ in $L(R^\vee)/L(R^\vee_S)$ is a torsion
element, say of order $d$. Any prime $p$ dividing $d$ is a torsion prime
for the root system $R$. By Corollary 2.8 of \cite{steinberg75} $p$
divides $|W|$, and therefore $p$ is invertible in $k$. We
conclude that $d$ is invertible in $k$. 

Now $d\alpha^\vee \in L(R^\vee_S)$, and therefore 
\[
dd_\alpha\alpha=\varphi(d\alpha^\vee) \in L(R_S),
\]
which implies that 
 $dd_\alpha\alpha(u)=0$ for all $u \in S$. Since
$dd_\alpha$ is invertible in $k$, we conclude that $\alpha(u)=0$ for all $u
\in S$, so that
$\alpha \in R_S$, as desired.

(3) It follows easily from \cite[Ch.~VI, no.~1.7, Prop.~24]{bourbaki.root}
that the $\mathbb Q$-closed subsets of $R$ are precisely those of the form
$R_M$ for some Levi subgroup $M \supset T$. 
\end{proof}

\subsection{A property of $\mathfrak a_M$} Let $M$ be a Levi subgroup of
$G$ containing $T$. Let $R_M$ be the set of roots of $T$ in $M$. Define a
linear subspace $\mathfrak a_M$ of $\mathfrak t$ by 
\[
\mathfrak a_M:=\{ u \in \mathfrak t: \alpha(u)=0 \quad \forall \, \alpha
\in R_M \}. 
\]

\begin{lemma}\label{lem.goth.am}
Assume, as usual, that $|W|$ be invertible in $k$. Then 
\[
R_M=\{\alpha \in R: \alpha(u)=0 \quad \forall \, u \in \mathfrak a_M \}. 
\]
\end{lemma}
\begin{proof}
Choose a base $B$ for the root system $R$ in such a way that $B \cap R_M$
is a base for $R_M$. Since the index of connection of $R$ divides $|W|$
(see \cite[Ch.~VI, no.~2.4, Prop.~7]{bourbaki.root}) and is therefore
invertible in
$k$, the elements in
$B$ yield linearly independent elements of $\mathfrak t^*$. Note that
$\mathfrak  a_M$ is the intersection of the root hyperplanes in $\mathfrak
t$ determined by the elements in $B\cap R_M$. 

We must show that if $\alpha \in R \setminus R_M$, then $\alpha $ does not
vanish identically on $\mathfrak a_M$. We may assume that $\alpha$ is
positive. Inside $X^*(T)$ we write $\alpha$ as a $\mathbb Z$-linear
combination of elements in $B$. Then some element $\beta \in B,\beta \notin
R_M$ occurs in this linear combination with positive coefficient $n$. It is
enough to show that
$n$ is non-zero in $k$. This is clear unless $k$ has characteristic $p$ for
some prime $p$.

Let $n'$ be the coefficient of $\beta$ in the highest root $\tilde\alpha$.
Then $n \le n'$, so it is enough to show that $n'<p$. This follows from our
hypothesis that $|W|$ be invertible in $k$ (check case-by-case).
\end{proof}

\section{Appendix: Fixed points of the action of a finite group on a
scheme}\label{sec.fp.gx}  
Throughout this section $G$ denotes a finite group of order  
$|G|$. For any set $Z$ on which $G$ acts we write
$Z^G$ for the set of fixed points of the action of $G$ on $Z$. Finally, $S$
denotes some scheme, which will often serve as a base scheme.  

\subsection{Review of coinvariants of $G$-actions on quasicoherent sheaves}
Let $X$ be a scheme and $\mathcal F$ a quasicoherent $\mathcal
O_X$-module.  We consider an action of $G$ on $\mathcal F$, in other words,
a homomorphism 
$\rho: G \to \Aut_{\ox}(\mathcal F)$. 

We write $\mathcal F_G$ for the \emph{coinvariants} of $G$ on $\mathcal F$.
By definition $\mathcal F_G$ is the quasicoherent $\ox$-module obtained as
the cokernel of the homomorphism 
\[
\bigoplus_{g \in G}\mathcal F \to \mathcal F
\]
whose restriction to the summand indexed by $g \in G$ is
$\rho(g)-\id_\mathcal F$. 

For any $\ox$-module $\mathcal H$ there is an obvious action of $G$ on
$\Hom_{\ox}(\mathcal F,\mathcal H)$, and it is evident from the definition
of coinvariants that there is a canonical isomorphism 
\begin{equation}\label{eq.coinv}
\Hom_{\ox}(\mathcal F_G,\mathcal H)=\bigl(\Hom_{\ox}(\mathcal F,\mathcal H)
\bigr)^G.
\end{equation}

\subsection{Fixed points of $G$-actions on schemes} 
Let $X$ be a scheme over $S$. Suppose that the finite group
$G$ acts on $X$ over $S$, by which we mean that for each $g \in G$ the
morphism $x \mapsto gx$ from $X$ to itself is a morphism over $S$. 

We define a contravariant set-valued functor $X^G$ on the category of
schemes $T$ over $S$  by the rule 
\[
X^G(T):=X(T)^G.
\]
\begin{lemma}\label{lem.diag}\hfill
\begin{enumerate}
\item The subfunctor $X^G$ of $X$ is represented by a locally closed
subscheme of $X$.
\item If $X$ is separated over $S$, then $i:X^G \hookrightarrow X$ is a
closed immersion.
\item If $X$ is locally of finite presentation over $S$, then $X^G$ is
locally of finite presentation over $S$. 
\item Taking fixed points commutes with arbitrary base change $S' \to S$,
which is to say that 
\[
(X\times_S S')^G=X^G \times_S S'.
\]
\end{enumerate}
\end{lemma}
\begin{proof}
Enumerate the elements of $G$ as $g_1,\dots,g_n$. Write $X^n$ for the
$n$-fold fiber product $X\times_S X \times_S \dots \times_S X$. We consider
two morphisms $X \to X^n$, one being the diagonal morphism $\Delta$ defined
by $\Delta(x)=(x,\dots,x)$, the other, denoted $\alpha$, being defined by
$\alpha(x)=(g_1x,\dots,g_nx)$. Taking the fiber product of these two
morphisms, we get a scheme over $S$ which clearly represents $X_G$. 

Thus we have a cartesian square
\[
\begin{CD}
X^G @>i>> X \\
@VVV @VV\alpha V \\
X @>\Delta>> X^n
\end{CD} 
\] 
showing that any property of $\Delta$ which is stable under base change
will be inherited by $i$. This proves (1), (2) and reduces (3) to
checking that $\Delta$ is locally of finite presentation when $X$ is
locally of finite presentation over $S$. This follows immediately from EGA
IV (1.4.3)(v), applied to the composition $\mathrm{pr}_1 \circ \Delta$,
with 
$\mathrm{pr}_1:X^n \to X$ denoting projection on the first factor. 

Finally, (4) is obvious from the definition of $X^G$. 
\end{proof}

\subsection{$1$-forms over fixed point subschemes}
For any scheme $X$ over $S$ one has the quasicoherent
$\mathcal O_X$-module $\Omega^1_{X/S}$ of $1$-forms on $X/S$, as well as
the tangent ``bundle'' $T_{X/S}$, which is the scheme, affine over $X$,
obtained as the spectrum of the symmetric algebra on the $\mathcal
O_X$-module $\Omega^1_{X/S}$. 

Consider a morphism $f:Y \to X$ of schemes over $S$, and a 
quasicoherent $\mathcal O_Y$-module $\mathcal H$. We regard $\mathcal O_Y
\oplus \mathcal H$ as an $\mathcal O_Y$-algebra in the usual way:
\[
(a_1,h_1)\cdot(a_2,h_2)=(a_1a_2,a_1h_2+a_2h_1).
\]
Put $Y(\mathcal H):=\Spec(\mathcal O_Y \oplus \mathcal H)$, a scheme affine
over $Y$. The augmentation $\mathcal O_Y \oplus \mathcal H \to \mathcal
O_Y$ (sending $(a,h)$ to $a$) yields a section of $Y(\mathcal H) \to Y$,
which we use to identify
$Y$ with a closed subscheme of $Y(\mathcal H)$ having the same underlying
topological space as $Y(\mathcal H)$. We then have (see EGA IV, 16.5) the
following  property of $\Omega^1_{X/S}$:
\begin{equation}\label{eq.om.ump}
\Hom_{\mathcal O_Y}(f^*\Omega^1_{X/S},\mathcal H)=\{\tilde f \in
\Hom_S(Y(\mathcal H),X):\tilde f|_Y=f\}.
\end{equation}

\begin{lemma}
Let $X$ be a scheme over $S$, and suppose that the finite group $G$ acts on
$X$ over $S$. Let $i:X^G \hookrightarrow X$ be the obvious inclusion. Then
there are canonical isomorphisms
\begin{equation}\label{eq.omegaG}
\Omega^1_{X^G/S}=(i^*\Omega^1_{X/S})_G
\end{equation}
and 
\begin{equation}\label{eq.tanG}
T_{X^G/S}=(T_{X/S})^G.
\end{equation}
The subscript $G$ on the right side of \eqref{eq.omegaG} indicates that we
take coinvariants for the action of $G$. 
\end{lemma}
\begin{proof}
To prove \eqref{eq.omegaG} it is enough to construct, for any quasicoherent
$\mathcal O_{X^G}$-module $\mathcal H$, a functorial isomorphism
\[
\Hom_{\mathcal O_{X^G}}(\Omega^1_{X^G/S},\mathcal H)=\Hom_{\mathcal
O_{X^G}}((i^*\Omega^1_{X/S})_G,\mathcal H).
\]
By \eqref{eq.om.ump} we have 
\[
\Hom_{\mathcal
O_{X^G}}(i^*\Omega^1_{X/S},\mathcal H)=\{\tilde i \in
\Hom_S(X^G(\mathcal H),X):\tilde i|_{X^G}=i\}.
\]
Taking invariants under $G$ and using \eqref{eq.coinv}, we see that 
\[
\Hom_{\mathcal
O_{X^G}}((i^*\Omega^1_{X/S})_G,\mathcal H)=\{\tilde i \in
\Hom_S(X^G(\mathcal H),X^G):\tilde i|_{X^G}=\id_{X^G}\},
\]
and by \eqref{eq.om.ump} the right side of this equality is equal
to 
\[
\Hom_{\mathcal O_{X^G}}(\Omega^1_{X^G/S},\mathcal H),
\]
as desired. 

From the definition of $T_{X/S}$ we have, for any scheme $S'$ over $S$, the
equality 
\[
T_{X/S}(S')=\{(f,\beta):f \in \Hom_S(S',X),\,\beta \in \Hom_{\mathcal
O_X}(\Omega^1_{X/S},f_*\mathcal O_{S'}) \}. 
\]
Taking fixed points under $G$, we find that 
\[
(T_{X/S})^G(S')=\{(f,\beta):f \in \Hom_S(S',X^G),\,\beta \in \Hom_{\mathcal
O_X}(\Omega^1_{X/S},i_*f_*\mathcal O_{S'})^G\}. 
\]
Using \eqref{eq.omegaG}, \eqref{eq.coinv} and the adjointness of $i^*$, 
$i_*$, we see that 
\[
\Hom_{\mathcal O_{X^G}}(\Omega^1_{X^G/S},f_*\mathcal O_{S'})=
\Hom_{\mathcal
O_X}(\Omega^1_{X/S},i_*f_*\mathcal O_{S'})^G, 
\]
from which it follows that 
\begin{align*}
(T_{X/S})^G(S')&=\{(f,\beta):f \in \Hom_S(S',X^G),\,\beta \in \Hom_{\mathcal
O_{X^G}}(\Omega^1_{X^G/S},f_*\mathcal O_{S'})\} \\
&=T_{X^G/S}(S'),
\end{align*}
which proves \eqref{eq.tanG}. 
\end{proof}

\subsection{Smoothness of fixed point subschemes}
Again consider an action of the
finite group $G$  on a scheme $X$ over $S$.
\begin{lemma}
Suppose that $X$ is smooth over $S$ and that $|G|$ is invertible on $S$.
Then $X^G$ is smooth over $S$.
\end{lemma}
\begin{proof}
It follows from Lemma \ref{lem.diag} (3) that $X^G$ is locally of finite
presentation over $S$. It remains to verify that $X^G$ is formally smooth
over $S$, so consider an affine scheme $\Spec(A)$ over $S$ and an ideal $I
\subset A$ such that $I^2=0$. Writing $X^G(A)$ for $\Hom_S(\Spec A,X^G)$,
we must show that
\[
\alpha:X^G(A) \to X^G(A/I)
\] 
is surjective.

Since $X$ is smooth over $S$, we do know that
\[
\beta:X(A) \to X(A/I)
\] 
is surjective. Given $x \in X(A/I)$, in other words an $S$-morphism $x:\Spec
A/I \to X$, the fiber $\beta^{-1}(x)$ is a principal homogeneous space
under (again see EGA IV, 16.5)
\[
M:=\Hom_{A/I}(x^*\Omega^1_{X/S},I). 
\]
Now suppose that $x \in X^G(A/I)$. Then $G$ acts compatibly on $M$ and
$\beta^{-1}(x)$, and the obstruction to the existence of a $G$-invariant
element in $\beta^{-1}(x)$ lies in $H^1(G,M)$. Since $M$ is a $G$-module on
which multiplication by $|G|$ is invertible, all higher group cohomology of
$M$ vanishes, so our obstruction is automatically trivial. Therefore 
$\alpha^{-1}(x)$ is non-empty, showing that $\alpha$ is surjective, as
desired. 
\end{proof}

In the next result we no longer need a base scheme $S$. Note that any
action of a finite group on a scheme $X$ is automatically an action on $X$
over $\Spec(\mathbb Z)$, so $X^G$ still makes sense and is a scheme (over
$\Spec(\mathbb Z)$). 
\begin{lemma}\label{lem.Gsmth}
Suppose that the finite group $G$ acts on schemes $X$, $Y$. Suppose further
that we are given  a $G$-equivariant morphism $f:Y \to X$, and consider
the morphism $Y^G \to X^G$ induced by $f$. 
\begin{enumerate}
\item If $Y$ is locally of finite presentation over $X$, then $Y^G$ is
locally of finite presentation over $X^G$. 
\item There is a canonical isomorphism 
$\Omega^1_{{Y^G/X^G}}=(i^*\Omega^1_{Y/X})_G$, where $i$ denotes the
inclusion $Y^G \hookrightarrow Y$ and the subscript
$G$ indicates coinvariants.  
\item If $Y$ is smooth over $X$, and $|G|$ is invertible on $X$, then $Y^G$
is smooth over $X^G$. 
\item If $Y \to X$ is a locally closed immersion, then so is $Y^G \to X^G$.
\end{enumerate} 
\end{lemma}
\begin{proof}
 We already know the first three parts of the lemma when $G$ acts
trivially on
$X$, so that $X^G=X$. To treat the general case, we form the cartesian
square 
\[
\begin{CD}
Y' @>>> Y \\
@VVf'V @VVfV \\
X^G @>>> X.
\end{CD} 
\]
The group $G$ still acts on the locally closed subscheme $Y'$ of $Y$, and it
is clear that $(Y')^G=Y^G$. If $f$ is locally of finite presentation
(respectively, smooth), then $f'$ is locally of finite presentation
(respectively, smooth). Moreover
$i^*\Omega^1_{Y/X}=(i')^*\Omega^1_{Y'/X^G}$, where $i'$ is the inclusion
$Y^G \hookrightarrow Y'$. Therefore it is enough to prove  the first
three parts of the lemma with
$f$ replaced by $f'$, and then we are done by the remark made at the
beginning of the proof. 

We now prove the last part of the lemma. Using that $f$ is a monomorphism,
we see that the square 
\[
\begin{CD}
Y^G @>>> Y \\
@VVV @VVfV \\
X^G @>>> X.
\end{CD} 
\]
is cartesian, allowing us to deduce that $Y^G \to X^G$ is an immersion from
the fact that $f$ is an immersion. 
\end{proof}

\section{Appendix. Greenberg's functor}\label{app.gre}
\subsection{Definition of Greenberg's functor}\label{sub.gre.sm} 
Let $X$ be a scheme of finite type over $\mathcal O$, and let $N$ be a
positive integer. Then Greenberg's functor associates to $X$ the  
scheme $X_N$ of finite type over $k$ whose
points in any $k$-algebra $A$ are given by 
\begin{equation}\label{eq.def.gre}
X_N(A):=X(A\otimes_k(\mathcal
O/\epsilon^N\mathcal O)).
\end{equation}
In particular the set of $k$-points of $X_N$ 
is $X(\mathcal O/\epsilon^N\mathcal O)$. 

An $\mathcal O$-morphism $f:Y \to X$ between schemes of finite type
over $\mathcal O$ induces a $k$-morphism 
\[
f_N:Y_N \to X_N.
\]

If $f$ is smooth (respectively,
\'etale), then
$f_N:Y_N \to X_N$ is smooth (respectively,
\'etale). Indeed, due to \eqref{eq.def.gre},  the formal
smoothness (respectively,
\'etaleness) of
$f_N$ is inherited from the formal smoothness (respectively,
\'etaleness) of $f$.

\subsection{The smooth case}\label{sub.gre.tan}
Suppose that $X$ is smooth over $\mathcal O$. Then $X_N$ is smooth over $k$.
It follows from
\eqref{eq.def.gre} and \eqref{eq.om.ump} that the tangent space to
$X_N$ at $x \in X_N(k)=X(\mathcal
O/\epsilon^N\mathcal O)$ is given by 
\[
T_{X_N,x} = x^*\mathcal T_{X/\mathcal O},
\]
where $x$ is being regarded as an $\mathcal O$-morphism $\Spec(\mathcal
O/\epsilon^N\mathcal O) \to X$, and $\mathcal T_{X/\mathcal O}$ is the
relative tangent sheaf of $X/\mathcal O$. Note that the tangent space
$T_{X_N,x}$ is in a natural way an $\mathcal O/\epsilon^N\mathcal
O$-module, free of finite rank. 

Suppose that $x$ is obtained by reduction modulo $\epsilon^N$ from $\tilde
x \in X(\mathcal O)$. Then $\tilde x^*\mathcal T_{X/\mathcal O}$ is a free
$\mathcal O$-module of finite rank that we will refer to informally as the
tangent space to $X(\mathcal O)$ at $\tilde x$ and denote by $T_{X(\mathcal
O),\tilde x}$.  Clearly we have 
\begin{equation}\label{eq.tan.iso}
T_{X_N,x}=T_{X(\mathcal
O),\tilde x}\otimes_\mathcal O (\mathcal O/\epsilon^N\mathcal O).
\end{equation}

Now suppose that $f:Y \to X$ is an $\mathcal O$-morphism between two smooth
schemes over $\mathcal O$; applying Greenberg's functor to $f$  we get a
$k$-morphism 
\[
f_N:Y_N \to X_N.
\]
 Let
$\tilde y
\in Y(\mathcal O)$ and put
$\tilde x:=f(y) \in X(\mathcal O)$; then let $y \in Y(\mathcal
O/\epsilon^N\mathcal O)$, $x \in X(\mathcal O/\epsilon^N\mathcal O)$ be the
points obtained from $\tilde y$, $\tilde x$ by reduction modulo
$\epsilon^N$. The differential of $f$ gives us an $\mathcal O$-linear map 
\[
df_{\tilde y}:T_{Y(\mathcal O),\tilde y} \to 
T_{X(\mathcal O),\tilde x}.
\]
Reducing this map modulo $\epsilon^N$ and using the isomorphism
\eqref{eq.tan.iso}, we obtain a $k$-linear map 
\[T_{Y_N, y} \to 
T_{X_N, x}
\]
which is easily seen to coincide with the differential of $f_N$ at $y$. 
In
other words, the differential of $f_N$ is the reduction modulo $\epsilon^N$
of the differential of $f$. 

\subsection{Restriction of scalars $\mathcal O_E/\mathcal O$ and
Greenberg's functor} \label{sub.gre.re.sc} 
Let $E$ be a finite extension  field of $F$, and let
$\mathcal O_E$ be the integral closure of $\mathcal O$ in $E$. 

Let $X$ be a scheme of finite type over $\mathcal O_E$. We denote by
$R_{\mathcal O_E/\mathcal O} X$ the scheme of finite type over $\mathcal O$
obtained by (Weil) restriction of scalars from $\mathcal O_E$ to $\mathcal
O$. Recall that the points of  
$R_{\mathcal O_E/\mathcal O} X$ in any $\mathcal O$-algebra $A$ are given
by 
\[
(R_{\mathcal O_E/\mathcal O} X)(A)=X(A \otimes_\mathcal O \mathcal O_E).
\]

Let $N$ be a positive integer. Applying Greenberg's functor to $R_{\mathcal
O_E/\mathcal O} X$ provides us with a
$k$-scheme whose set of $k$-points is $(R_{\mathcal O_E/\mathcal O}
X)(\mathcal O/\epsilon^N \mathcal O)=X(\mathcal O_E/\epsilon^N \mathcal
O_E)$. 
But there is another equally natural way to produce a $k$-scheme with the
same set of $k$-points, namely to apply Greenberg's functor (for the field
$E$ rather than the field $F$)  to
$X$ (and the quotient ring $\mathcal O_E/\epsilon^N \mathcal O_E$ of
$\mathcal O_E$). 
In fact these two $k$-schemes are canonically isomorphic, since for both
schemes the set of $A$-valued points ($A$ now being a $k$-algebra) works out
to be 
\[
X(A \otimes_k(\mathcal O_E/\epsilon^N \mathcal O_E)).
\]

\subsection{Greenberg's functor and fixed point sets} \label{sub.gre.fps}
Let $X$ be a scheme
of finite type over $\mathcal O$, and $N$ a positive integer. From
Greenberg's functor we get the $k$-scheme $X_N$.  Now suppose 
further that we are given an action of a finite group
$G$ on $X$ over
$\mathcal O$. Then, by functoriality, $G$ acts on $X_N$ over $k$. 

It follows immediately from the definitions that $(X_N)^G$ is canonically
isomorphic to $(X^G)_N$. Indeed, for both schemes the set of $A$-valued
points ($A$  being a $k$-algebra) works out to be 
\[
X(A \otimes_k(\mathcal O/\epsilon^N \mathcal O))^G.
\]

\bibliographystyle{amsalpha}

%\bibliography{thebib}

\providecommand{\bysame}{\leavevmode\hbox to3em{\hrulefill}\thinspace}

\end{document}